\documentclass[journal]{IEEEtran}

\usepackage{xspace,amsmath,amssymb,amsthm,amsfonts,epsfig,syntonly,dsfont,pifont}
\usepackage{cite,bm,color,url,textcomp}
\usepackage{algorithmicx,algorithm,algpseudocode}
\usepackage{epstopdf}
\usepackage{empheq,float}
\usepackage{graphicx,subfigure,balance}
\usepackage{multirow, multicol,booktabs}

\usepackage{makecell}

\usepackage[dvipsnames]{xcolor}

\newtheorem{assumption}{Assumption}
\newtheorem{assumptionalt}{Assumption}[assumption]

\newenvironment{assumptionp}[1]{
  \renewcommand\theassumptionalt{#1}
  \assumptionalt
}{\endassumptionalt}

\newtheorem{theorem}{Theorem}
\newtheorem{theoremalt}{Theorem}[theorem]

\newenvironment{theoremp}[1]{
  \renewcommand\thetheoremalt{#1}
  \theoremalt
}{\endtheoremalt}

\newtheorem{lemma}{Lemma}
\newtheorem{corollary}{Corollary}
\newtheorem{definition}{Definition}
\theoremstyle{definition}

\DeclareMathOperator*{\argmin}{arg\,min}
\DeclareMathOperator*{\tr}{tr}

\graphicspath{{./figs/}}

\usepackage{times}

\title{Conic Descent Redux \\for Memory-Efficient Optimization}

\author{Bingcong Li, and Georgios B. Giannakis\\

\thanks{
B. Li and G. B. Giannakis are with the Dept. of Electrical and Computer Engineering, University of Minnesota, Minneapolis, MN 55455 USA. Emails:~\{lixx5599, georgios\}@umn.edu.
}
}

\begin{document}

\maketitle
\begin{abstract}
Conic programming has well-documented merits in a gamut of signal processing and machine learning tasks. This contribution revisits a recently developed first-order conic descent (CD) solver, and advances it in three aspects: intuition, theory, and algorithmic implementation. It is found that CD can afford an intuitive geometric derivation that originates from the dual problem. This opens the door to novel algorithmic designs, with a momentum variant of CD, \underline{mo}mentum \underline{co}nic descent (MOCO) exemplified. Diving deeper into the dual behavior CD and MOCO reveals: i) an analytically justified stopping criterion; and, ii) the potential to design preconditioners to speed up dual convergence. Lastly, to scale semidefinite programming (SDP) especially for low-rank solutions, a memory efficient MOCO variant is developed and numerically validated.
\end{abstract}

\section{Introduction}\label{sec.intro}

Consider a conic programming setup of the form
\begin{align}\label{eq.prob}
	\min_{\mathbf{x} \in \mathbb{R}^d}  f(\mathbf{x}) ~~~\text{s.t.} ~ \mathbf{x} \in {\cal K} \end{align}
where the differentiable objective function $f$ is convex, and ${\cal K}$ denotes a convex cone. Conic problems are frequently encountered in machine learning and signal processing, where cones can for instance capture non-negative orthant constraints, second-order cones, positive semidefinite cones, exponential cones, and copositive cones \cite{mosek, boyd2004, dur2010}. The generality of conic problems fertilizes a number of application domains, leading to the well-documented success in applications such as community detection, and multi-task learning \cite{hajek2016,kato2007,hanasusanto2018conic}.

This work considers first order methods for solving \eqref{eq.prob}. We will focus on Frank Wolfe (FW) variants \cite{frank1956,jaggi2013,li2021heavy} since their computationally lightweight subproblems can avoid projection onto cones. Taking positive semidefinite cones $\mathbb{S}_+^n$ as an example, projection requires a full SVD with complexity ${\cal O}(n^3)$, while a FW subproblem only needs to find out the eigenvector associated with largest eigenvalue for certain matrix, reducing the overall complexity to ${\cal O}(n^2)$.

Nonetheless, the noncompact cone constraint prevents applying FW directly on problem \eqref{eq.prob}. A straightforward approach is to include a manually designed constraint to shrink the original constraint set to a compact one $\tilde{\cal K}$. Consider a simple example with ${\cal K} = \{ (x, y)| x \geq 0, y\geq 0\}$, one manner to define $\tilde{\cal K}$ is to turn the non-negative orthant into a polyhedron by adding another constraint, e.g., $x + y \leq 1$. This idea is formalized and generalized in \cite{harchaoui2015}, yet prior knowledge is of critical importance to the shrunk constraint $\tilde{\cal K}$ otherwise it may not contain optimal solutions to \eqref{eq.prob}.

Work \cite{locatello2017greedy} considers problem \eqref{eq.prob} with ${\cal K}$ having a relatively simple atomic expression. This additional assumption on ${\cal K}$ wipes out constraints such as doubly nonnegative cones, which are useful for reformulating combinatorial problems \cite{cui2020projecting}. Moreover, the convergence rate in \cite{locatello2017greedy} depends on the geometry of the cone, thus can be challenged by some ``illy conditioned'' ones. 

A recent method \cite{duchi2020conic} develops conic descent (CD) that can cope with general convex cones regardless of the atomic form of ${\cal K}$. However, many of first order approaches, including CD, only target at the primal convergence, leaving the dual properties relatively untouched despite conic duality can be informative. In this work, a detailed study is carried out to understand the dual convergence of CD. In particular, we first provide an explanation to CD that is not only geometrically intuitive, but also having matching mathematical support in the dual domain. This explanation brings up opportunities on algorithmic design, and resulting in a new variant of CD, \underline{mo}mentum \underline{co}nic descent (MOCO). MOCO is equipped with heavy ball momentum for faster convergence. Then, extensive theoretical analyses on dual domain bring up deeper insights, and a practical stopping criterion to estimate suboptimality.

We then focus on an instance of \eqref{eq.prob}, SDP problems \cite{vandenberghe1996, nesterov2000semidefinite} with the goal of improved scalability. The key observation that motivates the study of memory efficient SDP is that many SDP instances are raised up from vector problems. We term such problems as \textit{raised-up SDPs.} Consider a simple quadratic problem $\min_{\mathbf{x} \in \mathbb{R}^d} \| \mathbf{x} \|_2^2$ as an example. Upon letting $\mathbf{X}:= \mathbf{x}\mathbf{x}^\top$, one can rewrite this problem as $\min_{\mathbf{X}} \tr(\mathbf{X})~ \text{s.t.}~ \mathbf{X} \in \mathbb{S}_d^+$ and $\text{Rank}(\mathbf{X}) \leq 1$. Dropping the rank constraint, one ends up with a SDP problem. While the previous example is too simple to visualize the benefit of raising up vector problem to SDPs, often times such a technique is helpful to turn a nonconvex problem into a convex one; see e.g., \cite{wang2015robust,yu2004iterative} for more benefits in real-world applications. However, the raised-up SDP is at an obvious cost of increasing storage relative to its vector form. Our goal is to alleviate such a memory issue leveraging the observation that the desirable solution is usually low rank (recall the rank constraint in our toy example). We propose a memory efficient implementation of MOCO, and leverage Burer-Monteiro (BM) factorization heuristic \cite{burer2003} to further enhance its empirical performances.

 In succinct form, our contributions are listed as follows.

\begin{enumerate}
	\item[\ding{118}] It is found that conic descent (CD) admits a geometrical explanation. Interestingly, the geometry has a rigorous mathematical foundation in the dual domain of \eqref{eq.prob}. 
	
	\item[\ding{118}] A new algorithm is developed based on the geometrical interpretation. The resultant approach, MOCO, improves the convergence rate of CD, and showcases numerical merits on tested problems.
	
	\item[\ding{118}] Comprehensive analyses to the dual properties are provided for MOCO. It is observed that the primal and dual convergence do not share the same rate, and the dual behavior can be influenced via preconditioning.
	
	\item[\ding{118}] We further modify MOCO for memory efficiency of large-scale (raised-up) SDPs. BM heuristic is also incorporated into modified MOCO to facilitate numerical performances. 
	
\end{enumerate}

\textbf{Notational conventions}. Bold lowercase (capital) letters denote column vectors (matrices); $\|  \cdot \|$ stands for a norm of either a vector or a matrix, whose dual norm is denoted by $\|\cdot \|_*$; and $\langle \cdot, \cdot \rangle$ is the inner product. Given a cone ${\cal K}$, its dual cone is written as ${\cal K}^*$. For a set ${\cal S}$, we let $\text{dist}(\mathbf{x}, {\cal S})$ and $\text{dist}_*(\mathbf{x}, {\cal S})$ denote the distance of vector $\mathbf{x}$ to set ${\cal S}$ w.r.t. $\|\cdot \|$, and $\| \cdot \|_*$, respectively. We use $\mathbb{S}^n$ for symmetric real matrices, and $\mathbb{S}^n_+$ to denote the semidefinite positive cone, i.e., all symmetric real positive semidefinite matrices of size $n \times n$.

\section{Understanding Conic Descent Geometrically}

This section first describes in detail the class of problems that we are interested in, and then exemplifies a $2$-dimensional toy problem to unveil the underlying intuition of CD.

\subsection{Basic assumptions}\label{sec.assumptions}

We formally pinpoint problem \eqref{eq.prob} by mildly confining the class of objective functions.

\begin{assumption}[Lipschitz continuous gradient]\label{as.1}
	The objective function $f: {\cal K}\rightarrow \mathbb{R}$ has $L$-Lipchitz continuous gradients; i.e., $\|\nabla f(\mathbf{x}) - \nabla f(\mathbf{y}) \|_* \leq L \| \mathbf{x}-\mathbf{y} \|, \forall\, \mathbf{x}, \mathbf{y} \in {\cal K}$.
\end{assumption} 

\begin{assumption}[Strictly convex loss]\label{as.2}
	The objective function $f: {\cal K} \rightarrow \mathbb{R}$ is strictly convex; that is, $f(\mathbf{y}) - f(\mathbf{x}) > \langle \nabla f(\mathbf{x}), \mathbf{y} - \mathbf{x} \rangle, \forall\, \mathbf{x}, \mathbf{y} \in {\cal X}$ where $\mathbf{x} \neq \mathbf{y}$.
\end{assumption} 

Assumption \ref{as.1} is standard in optimization literatures \cite{jaggi2013,duchi2020conic,nesterov2004,li2020,li2021heavy,li2020extra,zhang2021}. Assumption \ref{as.2} is slightly stronger than the commonly adopted one that only requires convexity. This is because of the need of a regularity condition on $f$ to ensure the existence of an optimal solution. Although not stated, other works such as \cite{locatello2017greedy} also need this regularity conditions. For example, it is impossible to minimize $f(x,y) = -x + y^2$, which is not strictly convex, over the cone ${\cal K}:= \{ (x, y)|x \geq 0, y \geq 0 \}$. Nonetheless, Assumption \ref{as.2} is easily satisfied in practice, since it covers many prevalent loss functions, for example, squared $\ell_2$ loss and logistic loss. Note that Assumption \ref{as.2} is slightly stringent for SDPs, and we will relax it for a large class of SDPs later in Section \ref{sec.sdp}. It is also possible to regulate $f$ with assumptions other than strictly convex. For example, the work \cite{duchi2020conic} assumes $f$ to have no nonzero direction of recession in ${\cal K}$. Despite this assumption is difficult to verify in practice, our results extends to this setting after justifying the notation accordingly.

For the constraint, we also require the cone to be convex, implying convexity of \eqref{eq.prob}.
\begin{assumption}[Convex cone]\label{as.3}
	The constraint set ${\cal K} \in \mathbb{R}^d$ is a convex cone; i.e., $\lambda_1\mathbf{x} + \lambda_2 \mathbf{y} \in {\cal K}$ for any $\lambda_1 \geq 0$, $\lambda_2 \geq 0$ and $\mathbf{x}, \mathbf{y} \in {\cal K}$.
\end{assumption}

There are several natural approaches to solve \eqref{eq.prob} under Assumptions \ref{as.1} -- \ref{as.3}.

\textbf{Approach 1.} Projected gradient descent (GD) is perhaps the first idea coming into mind. The issue with GD, however, is that projection on a cone can be expensive; see the earlier example of semidefinite positive cone in Section \ref{sec.intro}.

\textbf{Approach 2.} If one has a hint of $\|  \mathbf{x}^*\|$, where $\mathbf{x}^*$ is an optimal solution to \eqref{eq.prob}, it is possible to manually impose compactness by including an additional constraint $\| \mathbf{x} \| \leq R$ to \eqref{eq.prob} to clear the obstacles of applying FW. While the FW subproblem is typically much cheaper than projection, a proper estimation on $\|  \mathbf{x}^*\|$ is challenging if not impossible. An overestimated $\|  \mathbf{x}^*\|$ degrades the performance of FW since its convergence is shaped heavily by the diameter of the constraint \cite{jaggi2013}; while an underestimated $\|  \mathbf{x}^*\|$ may exclude the optimal solution from the feasible domain.

Given the downside of these two approaches, there is a pressing need of more efficient methods. A recent work \cite{duchi2020conic} introduces conic descent (CD). However, the lack of intuition somehow shades the popularity of this approach with great potential. Next, we unveil CD's underlying geometry.

\subsection{Geometric interpretation for CD}

\begin{figure}[t]
	\centering
	\includegraphics[width=.25\textwidth]{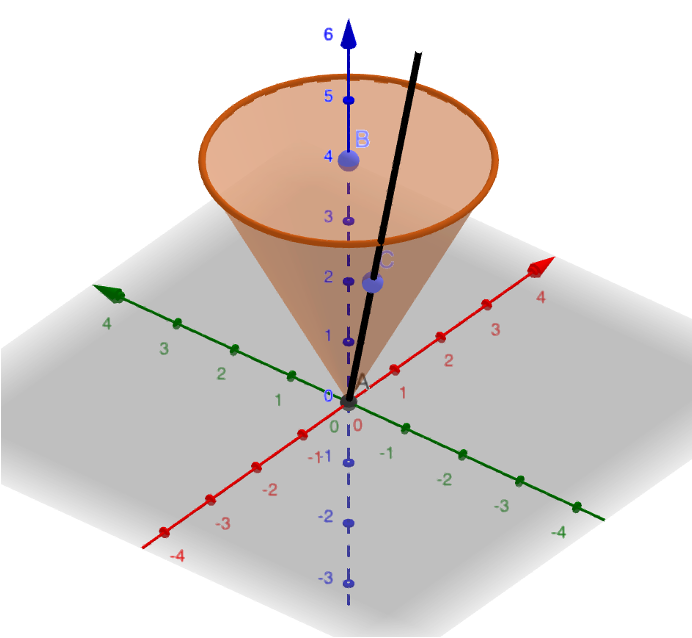}
	\caption{An example of an ice-cream cone.}
	 \label{fig.icecream}
\end{figure}

 Our novel interpretation of CD is built on two key observations. The first one is that a convex cone can be viewed as a set of rotated rays. We will only consider rays initialed at $\mathbf{0}$, that is, $\{t \mathbf{x}| t \geq 0\}$ for some $\mathbf{x} \neq \mathbf{0}$. For example, the first orthant in a $2d$-Cartesian plane can be viewed as the area scanned over by spinning the ray $\{(x, y)| x \geq 0, y = 0\}$ counterclockwise by $90$ degrees. Another example can be visualized in the ice-cream cone in Figure \ref{fig.icecream}. The second observation is that minimizing over a ray is an $1d$-convex problem and can be solved easily or even analytically. 
 
 Indeed, finding an optimal solution $\mathbf{x}^*$ to \eqref{eq.prob} amounts to finding the ray containing it. While it is challenging to find the desirable ray in just a single step, one may progressively improve the \textit{quality} on a ray, which is defined as the minimum function value of this ray. This intuition prompts us to decouple \eqref{eq.prob} into two (series of) subproblems: i) ray search, where the goal is to guess a ray that may contain $\mathbf{x}^*$; and ii) ray minimization, where this ray is minimized to obtain its quality. The overall goal is that the quality of a ray is improved iteratively until the optimal ray is found. It turns out that CD follows exactly this iterative procedure. 

\textbf{A polar-coordinate perspective.} The previous intuition can be understood more concretely through a $2$-dimensional example. Consider a simple quadratic objective function 
	\begin{align*}
		f(x, y) = (x-1)^2 + y^2
	\end{align*}
	with the cone constraint being the positive orthant, i.e., 
	 \begin{align*}
	 	{\cal K}:= \{(x, y)| x \geq 0, y \geq 0 \}.
	 \end{align*}
This problem can be transformed into polar coordinate, where $(x, y)$ are substituted to angular variables $(r, \theta)$, where $r \in [0, +\infty)$, and $\theta \in [0, \frac{\pi}{2}]$. Defining $t:= \cos \theta$, we can further change the variables as $x = r t$ and  $y = r \sqrt{1-t^2}$. The problem therefore becomes
\begin{align}\label{eq.polar}
	\min_{r,t}~& g(r, t) := (rt - 1)^2 + r^2(1 - t^2) \\
	\text{s.t.}~& r \geq 0, ~t \in [0, 1]	 \nonumber
\end{align}
From this reformulation \eqref{eq.polar}, it is clear that ray search targets at the optimal $t^*$, and ray minimization is used to obtain $r^*$ on the previously found ray. 

\textbf{Issues of working on polar coordinate}. Despite the reformulation \eqref{eq.polar} is geometrically intuitive, challenges remain even for this toy example. The first difficulty comes from the fact that problem \eqref{eq.polar} is not necessarily convex as it is straightforward to verify that the Hessian of $g$ is negative-definite, i.e.,
\begin{align*}
	\nabla^2 g = \begin{bmatrix}
		2  & -2 \\
		-2 &  0 
	\end{bmatrix}.
\end{align*}
Secondly, it is not always easy to reformulate a problem to its polar form, especially for those high dimensional cases. Therefore, it is more attractive to work with non-reformulated form \eqref{eq.prob}, performing ray search in an implicitly manner through the key message from \eqref{eq.polar}, that is, \textit{ray search is essentially a problem on compact domain} ($t \in [0,1]$).

\subsection{The MOCO algorithm}

The conic problem \eqref{eq.prob} can be solved by alternating between ray search and ray minimization as explained in previous subsection. In contrast with CD that adopts vanilla FW for ray search \cite{duchi2020conic}, here we propose to augment ray search with with momentum FW \cite{li2021heavy}. The resultant approach, MOCO, is summarized in Alg. \ref{alg.mcd}. While ray minimizing is straightforward in line 3, ray search is more involved; see lines 4 -- 7. Note that MOCO boils down to CD in \cite{duchi2020conic} if $\delta_k \equiv 1$.

It is known that $\nabla f(\mathbf{x}_k)$ is not the best coefficient to use in FW subproblems \cite{li2021heavy,li2020extra}. This motivates the use of the heavy ball momentum in MOCO. MOCO subproblem in line 5 instead relies on $\mathbf{g}_k$, a weighted average of past gradients. The average $\mathbf{g}_k$ smoothes the possible rapid changes of gradients in consecutive iterations, leading to a more stable searching direction. Another benefit of using momentum is the possibility to continue optimizing even if $\theta_k=0$. This can be helpful for (matrix SDP) problems with structural solutions, e.g., sparsity or low rankness. The CD iteration stops if $\theta_k = 0$ (since later iterations will not move), and an optimal solution is thus found. On the other hand, the heavy ball momentum further adjusts the weight for $\nabla f(\mathbf{x}_k)$ and continues optimizing. To see this, note that when $\theta_k = 0$, we have $\eta_{k+1}\mathbf{x}_{k+1} = \eta_k \mathbf{x}_k$. As a result, $\mathbf{g}_{k+2} = (1-\delta_{k+1})\mathbf{g}_{k+1} + \delta_{k+1} \nabla f(\eta_{k+1}\mathbf{x}_{k+1}) = (1-\delta_{k+1})\mathbf{g}_{k+1} + \delta_{k+1} \nabla f(\eta_k\mathbf{x}_k)$, that is, the weight on $\nabla f(\eta_k \mathbf{x}_k)$ is adaptively increased to $\delta_k(1-\delta_{k+1}) + \delta_{k+1}$ if one further unpacks $\mathbf{g}_{k+1}$. This gives a different search direction to continue the search for e.g., lower rank solutions.

Different from standard FW subproblems, which is $\argmin_{\mathbf{v} \in {\cal K}} \langle \mathbf{g}_k, \mathbf{v} \rangle$, the MOCO subproblem (for ray search) adds an additional constraint $\| \mathbf{v} \| \leq 1$ to ensure that the subproblem is solvable. Concretely, this amounts to our constraint $t \in [0, 1]$ in the toy example \eqref{eq.polar}, and the additional constraint can be taken as the range on the cosine of angular variables. Adding the additional constraint $\| \mathbf{v} \| \leq 1$ typically induces no extra computational burden compared to FW subproblems. For example in the SDPs considered later, the subproblems of MOCO and FW have the same complexity. More on ray search will be discussed in Section \ref{sec.dual}, where we will view ray search from a duality lens.

\begin{algorithm}[t]
    \caption{\textbf{Mo}mentum \textbf{co}nic descent (MOCO)}\label{alg.mcd}
    \begin{algorithmic}[1]
    	\State \textbf{Initialize:} $\mathbf{x}_0 $, $\delta_k = \frac{2}{k+2} \forall k$
			\For {$k=0,1,\dots,K$}
				\State $\eta_k = \argmin_{\eta \geq 0} f(\eta \mathbf{x}_k)$ \Comment{Ray minimization}
				\State $\mathbf{g}_k = (1 - \delta_k) \mathbf{g}_{k-1} + \delta_k \nabla f( \eta_k \mathbf{x}_k) $
				\State $\mathbf{v}_k = \argmin_{\mathbf{v}} \langle \mathbf{g}_k, \mathbf{v} \rangle$ s.t. $\| \mathbf{v} \| \leq 1, \mathbf{v} \in {\cal K}$ 
				\State $\theta_k = \argmin_{\theta \geq 0} f(\eta_k\mathbf{x}_k + \theta \mathbf{v}_k)$ 
				\State $\mathbf{x}_{k+1} = \eta_k \mathbf{x}_k + \theta_k \mathbf{v}_k	$ \Comment{Ray search}
			\EndFor
		\State \textbf{Return:} $\eta_{K} \mathbf{x}_{K}$
	\end{algorithmic}
	
\end{algorithm}

\section{Primal-dual convergence of MOCO}

Having explained the intuition of MOCO, we next focus on its theoretical properties. It is not difficult to see that MOCO converges after the first iteration for any initialization $\mathbf{x}_0 \in {\cal K}$ if $\mathbf{x}^* = \mathbf{0}$. We will hence cope with nontrivial problems assuming $\mathbf{x}^* \neq \mathbf{0}$ in the following subsections.

\subsection{Primal convergence}

We first deliver a direct result of ray minimization.

\begin{lemma}\label{lemma.raymin}
	For every iteration, MOCO ensures that
	\begin{align}\label{eq.zero}
		\langle \eta_k \mathbf{x}_k, \nabla f(\eta_k \mathbf{x}_k) \rangle = 0.	
	\end{align}
\end{lemma}

Another preparation for the convergence proof is a series of helper functions $\Phi_{k+1}(\mathbf{x})$, defined as
\begin{align}\label{eq.phi}
	 \Phi_{k+1}(\mathbf{x}) := & (1-\delta_k) \Phi_k(\mathbf{x}) \\
	& + \delta_k \big[ f(\eta_k\mathbf{x}_k) + \langle \nabla f(\eta_k\mathbf{x}_k), \mathbf{x} \rangle \big], \forall k \geq 0. \nonumber
\end{align}
The definition of $\Phi_0(\mathbf{x})$ does not influence $\Phi_{k+1}(\mathbf{x})$ since $\delta_0 = 0$. Similar to the spirit of \cite{li2021heavy}, the helper functions can be regarded as lower bounds for $f(\mathbf{x})$, where the detailed implications can be found later in Lemma \ref{lemma.phi}. However, there is a key difference that brings up additional challenges to the analysis of MOCO. Unlike \cite{li2021heavy}, $\Phi_{k+1}(\mathbf{x})$ in \eqref{eq.phi} may have a minimum reaching $-\infty$ due to the noncompactness of ${\cal K}$. Consequently, one cannot adopt $\min_{\mathbf{x}\in {\cal K}} \Phi_{k+1}(\mathbf{x})$ directly as the lower bound for $f(\mathbf{x}^*)$. The remedy for this issue is summarized in the next lemma.

\begin{lemma}\label{lemma.phi}
	$\Phi_{k+1}(\mathbf{x})$ satisfies that: i) $\mathbf{v}_k$ minimizes $\Phi_{k+1}(\mathbf{x})$ over $\{\mathbf{x}| \mathbf{x} \in {\cal K}, \| \mathbf{x} \| \leq 1 \}$; and, ii) there exists $\rho_k \geq 0$ such that $f(\mathbf{x}^*) \geq \Phi_{k+1} ( \|\mathbf{x}^*\| \mathbf{v}_k) + \rho_k$ holds, where $\rho_k = 0$ only if $\{\eta_\tau \mathbf{x}_\tau \equiv \mathbf{x}^* \}_{\tau = 0}^k$. The rigorous expression of $\rho_k$ can be found in Appendix \ref{apdx.lemma.phi}.
\end{lemma}	

Lemma \ref{lemma.phi} shows that by concentrating on a region that is the intersection of ${\cal K}$ and a norm ball, minimizing $\Phi_k(\mathbf{x})$ enables an underestimate of $f(\mathbf{x}^*)$.		
		
\begin{theorem}[Primal convergence]\label{thm.primal}
Suppose that Assumptions \ref{as.1}, \ref{as.2}, and \ref{as.3} hold. Choosing $\delta_k = \frac{2}{k+2}$, MOCO in Alg. \ref{alg.mcd} guarantees that
	\begin{align*}
		f(\eta_{k+1} \mathbf{x}_{k+1}) - \Phi_{k+1}(\|\mathbf{x}^* \|\mathbf{v}_k)  \leq \frac{2L \|\mathbf{x}^* \|^2}{k+2} - \rho_k
	\end{align*}
	where $\rho_k$ is defined in Lemma \ref{lemma.phi}.
\end{theorem}

The convergence rate of MOCO can be established as a simple combination of Theorem \ref{thm.primal} and Lemma \ref{lemma.phi}.

\begin{corollary}
	Under assumptions and parameter choices in Theorem \ref{thm.primal}, Alg. \ref{alg.mcd} converges with a rate
	\begin{align*}
		f(\eta_{k+1}\mathbf{x}_{k+1}) 	- f(\mathbf{x}^*) \leq \frac{2L \|\mathbf{x}^* \|^2}{k+2} - \rho_k.
	\end{align*}
\end{corollary}

Comparing the rate of MOCO to its non-momentum counterpart, CD \cite{duchi2020conic}, it is observed that momentum tightens the convergence rate by a small term $\rho_k$. This validates the merits of applying momentum to ray search.

\subsection{Dual convergence}\label{sec.dual}

We then tackle the dual convergence of MOCO to gain a complete understanding of its behaviors. Note that our analysis techniques can be directly extended to CD \cite{duchi2020conic}, and each of theorem below has a CD counterpart that differs only in constant, which we omit to save space.

\begin{definition}
Let $\epsilon \geq 0$ be some desirable tolerance. A point $\mathbf{x}^*_\epsilon$ is said to satisfy KKT condition of \eqref{eq.prob} $\epsilon$-approximately if
\begin{subequations}\label{eq.kkt}
\begin{align}
	 & \mathbf{x}_\epsilon^* \in {\cal K} \label{eq.dual1} \\
	 & \langle \nabla f(\mathbf{x}_\epsilon^*), \mathbf{x}_\epsilon^* \rangle = 0 \label{eq.dual2}\\
	 &  \big[ \text{dist}_* \big(\nabla f(\mathbf{x}_\epsilon^*), {\cal K}^* \big) \big]^2 \leq \epsilon. \label{eq.dual3}
\end{align}
\end{subequations}
\end{definition}

In particular, \eqref{eq.dual1} denotes primal feasibility of $\mathbf{x}_\epsilon^*$, \eqref{eq.dual2} and \eqref{eq.dual3} characterize complementary slackness and dual feasibility, respectively. Note that the KKT condition is satisfied if $\big[ \text{dist}_* \big(\nabla f(\mathbf{x}_\epsilon^*), {\cal K}^* \big) \big]^2 = 0$ (i.e., $\epsilon = 0$). Our gaol here is to understand that how fast can $\{ \eta_k \mathbf{x}_k \}$ generated by MOCO converge to an $\epsilon$-approximate KKT point.

An obvious fact is that MOCO never generates points outside of ${\cal K}$. Hence, $\eta_k\mathbf{x}_k$ is always primal feasible with \eqref{eq.dual1} satisfied automatically. Equation \eqref{eq.dual2} is also satisfied by $\eta_k \mathbf{x}_k$ as a result of ray minimization; see Lemma \ref{lemma.raymin}. This further explains the role of ray minimization, that is, \textit{it seeks $\mathbf{x}^*$ by eliminating points that are not complementarily slack.} It turns out that $\{\eta_k \mathbf{x}_k\}$ is not always dual feasible. Hence, the number of iterations required to ensure dual feasibility characterizes how fast an $\epsilon$-approximate KKT solution is found. Toward this goal, the key inequality leveraged is summarized in the following lemma.

\begin{lemma}\label{lemma.smooth_implification}
	Suppose that Assumptions \ref{as.1} and \ref{as.2} hold, then we have
	\begin{align*}
		f(\mathbf{x}) - f(\mathbf{y}) \geq \langle \nabla f(\mathbf{y}), \mathbf{x} - \mathbf{y} \rangle + \frac{1}{2L} \|\nabla f(\mathbf{y}) - \nabla f(\mathbf{x})  \|_*^2.
	\end{align*}
\end{lemma}

Lemma \ref{lemma.smooth_implification} extends \cite[Theorem 2.1.5]{nesterov2004} to non-Euclidean norms, and it is critical to MOCO's dual convergence. 

\begin{theorem}[Dual convergence]\label{thm.dual}
	 Suppose that Assumptions \ref{as.1}, \ref{as.2} and \ref{as.3} hold. With $\delta_k = \frac{2}{k+2}$, MOCO guarantees that
	\begin{equation*}
		\big[ \text{dist}_* \big(\nabla f(\eta_k \mathbf{x}_k)), {\cal K}^* \big) \big]^2 \leq \frac{4L^2 \| \mathbf{x}^* \|^2}{k+1}.
	\end{equation*}
\end{theorem}

Theorem \ref{thm.dual} asserts that an $\epsilon$-approximate KKT point can be found by MOCO after at most ${\cal O}(\frac{L^2 \| \mathbf{x}^*\|^2}{\epsilon})$ iterations. A critical observation is that the $L$ dependence is different on primal (Theorem \ref{thm.primal}) and dual (Theorem \ref{thm.dual}). This difference will influence the design of stopping criterion, which will be discussed in detail in the upcoming subsection.

\subsection{Stopping criterion}

While Theorems \ref{thm.primal} and \ref{thm.dual} characterize the primal and dual convergence rates, it is still unclear that when is a good time to stop MOCO iteration. Simply setting $K = {\cal O}(\frac{1}{\epsilon})$ works, but it could be too pessimistic since the rates are established for worst cases. In this subsection, we pursue a quantifiable overestimate of suboptimality that not only converges to $0$ as $k$ grows, but also can be obtained as a byproduct of MOCO subproblem.

Stopping criterions can be designed based on either primal or dual errors. If working with the primal, $f(\eta_k \mathbf{x}_k) - \Phi_k(\|\mathbf{x}^* \|\mathbf{v}_{k-1})$ in Theorem \ref{thm.primal} can be leveraged as an optimality measure. However, its value is impossible to compute due to the lack of knowledge about $\|\mathbf{x}^*\|$. The attempt on dual domain is to rely on $\big[ \text{dist}_* \big(\nabla f(\eta_k\mathbf{x}_k), {\cal K}^* \big) \big]^2 $ in Theorem \ref{as.2}. The issue is, however, computationally expensive and impractical since it requires a projection onto ${\cal K}^*$. To overcome these limitations, we find that $\langle \mathbf{g}_k, \mathbf{v}_k \rangle$ approximates $\big[ \text{dist}_* \big(\nabla f(\eta_k\mathbf{x}_k), {\cal K}^* \big) \big]^2 $ well, and can be used as a tractable certification for optimality. 

To see this, we first write out the dual for MOCO subproblem in line 5,
\begin{align}\label{eq.line4dual}
	\max_{ \mathbf{u}} - \| \mathbf{g}_k - \mathbf{u} \|_* ~~ \text{s.t.}~~ \mathbf{u} \in {\cal K}^*.
\end{align}
This dual problem \eqref{eq.line4dual} projects $\mathbf{g}_k$ onto ${\cal K}^*$, and the optimal objective value is $-\text{dist}_*(\mathbf{g}_k, {\cal K}^*)$. Comparing with \eqref{eq.dual3}, it can be seen that long as $\text{dist}_*(\mathbf{g}_k, {\cal K}^*) \approx \text{dist}_*(\nabla f(\eta_k\mathbf{x}_k), {\cal K}^*)$, one can use the optimal value of \eqref{eq.line4dual} as stopping criterion. This observation is formalized in the following theorem.

\begin{theorem}[Stopping criterion]\label{thm.stopcr}
	Suppose that Assumptions \ref{as.1}, \ref{as.2}, and \ref{as.3} hold. Upon choosing $\delta_k = \frac{2}{k+2}$, the following inequality holds for MOCO in Alg. \ref{alg.mcd} for any $k \geq 2$
	\begin{align}\label{eq.tttt}
		\big{[}\text{dist}_*(\mathbf{g}_k, {\cal K}^*) \big{]}^2	\leq \frac{9.7  L^2\|  \mathbf{x}^*\|^2}{k+1}.
	\end{align}
\end{theorem}


Theorem \ref{thm.stopcr} shows that $\text{dist}_*(\mathbf{g}_k, {\cal K}^*)$ converges at the same rate of $\text{dist}_*(\nabla f(\eta_k\mathbf{x}_k), {\cal K}^*)$ up to constant factors. Hence it further gives a math interpretation for ray search, that is, \textit{it projects $\mathbf{g}_k$ to ${\cal K}^*$ for (approximated) dual feasibility}.

Next we show that $\text{dist}_*(\nabla f(\eta_k\mathbf{x}_k), {\cal K}^*)$ can be conveniently obtained, suiting for the need of the stopping criterion. Strong duality between \eqref{eq.line4dual} and line 5 means that $\text{dist}_*(\mathbf{g}_k, {\cal K}^*) = - \langle \mathbf{g}_k, \mathbf{v}_k \rangle$. Therefore, one can simply approximate $\text{dist}_*(\nabla f(\eta_k\mathbf{x}_k), {\cal K}^*)$ via $ \langle \mathbf{g}_k, \mathbf{v}_k \rangle$, and assert an $\epsilon$-approximated KKT point is found whenever
\begin{align}\label{eq.stop}
	\langle \mathbf{g}_k, \mathbf{v}_k \rangle \geq - {\cal O}(\sqrt{\epsilon}).
\end{align}

It worth pointing out that the criterion \eqref{eq.stop} is an estimation on dual feasibility, as oppose to the primal error $f(\mathbf{x}_k) - f(\mathbf{x}^*)$ in standard FW literatures \cite{jaggi2013,li2021heavy}. In other words, \eqref{eq.stop} is no longer affine invariant as in standard FW, opening the possibility for preconditioning.

\textbf{Preconditioning.} With the hope of faster numerical performance, preconditioning applies a linear transformation to $\mathbf{x}$ and solves the transformed problem. It is observed that preconditioning has different impacts on primal and dual of MOCO. In particular, precondition does not reduce $L\|\mathbf{x}^* \|^2$ in primal error (cf. Theorem \ref{thm.primal}), but it can shrink $L\|\mathbf{x}^* \|$ in dual error (cf. Theorems \ref{thm.dual} and \ref{thm.stopcr}). Consider the following simple example with $f(x) = (x-2)^2$, whose preconditioned version is given by $g(x) = f(2x) = (2x - 2)^2$. We will use subscript $f$ and $g$ to denote constants and variables related to $f(x)$ and $g(x)$, respectively. In this case, one can verify that $L_f \| x_f^* \|_2^2 = L_g \| x_g^* \|_2^2$, but $L_f \| x_f^* \|_2 \neq L_g \| x_g^* \|_2$, demonstrating that the dual error can be scaled down without affecting primal error via proper precondition schemes. Henceforth, an optimal preconditioner can reduce the value of stopping criterion, leading to faster termination of the iterative procedure. This gives new questions on how to find the best preconditioner, which we leave for future work.

\section{Memory efficient MOCO for SDPs}\label{sec.sdp}

To enhance practical merits, we further develop a specific implementation for MOCO to reduce the memory consumption of large-scale semidefinite programing (SDP). By leveraging the problem structure, it is possible not only to store vectors in lieu of full matrix variables, but also to relax the regularity condition, i.e., strict convexity, in Assumption \ref{as.2}. We also augment this memory efficient MOCO with a greedy step based on a Burer-Monteiro (BM) factorization heuristic. When injecting a greedy step, it usually improves MOCO convergence.

\subsection{Problem statement}

Consider SDPs of the following form
\begin{align}\label{eq.prob_mtrx}
	\min_{\mathbf{X} } ~ f\big({\cal G} ( \mathbf{X}) - \mathbf{z} \big) 
	 ~~~~\text{s.t.} ~ \mathbf{X} \in \mathbb{S}_+^n 
\end{align}
where $\mathbf{z} \in \mathbb{R}^d$ is a given vector. The linear operator ${\cal G}$ maps $\mathbf{X} \in \mathbb{S}_+^n$ to $\mathbb{R}^d$, and it is defined as
\begin{align}
	{\cal G} ( \mathbf{X}) := [\tr(\mathbf{G}_1 \mathbf{X}), \ldots, \tr(\mathbf{G}_d \mathbf{X})]^\top
\end{align}
where $\mathbf{G}_i \in \mathbb{S}^n, i = 1, \ldots, d$. Problem \eqref{eq.prob_mtrx} appears frequently in machine learning and statistics, where $\{\mathbf{G}_i\}$ are often structural, e.g., low rank, sparse, or discrete Fourier transformation matrices. Given ${\cal G}$, its adjoint on a vector $\mathbf{a} \in \mathbb{R}^d$ is 
\begin{align}
{\cal G}^*(\mathbf{a}) = a_1 \mathbf{G}_1 + \ldots + a_d \mathbf{G}_d .
\end{align}
In the sequel, we assume that $n^2 \gg d$, and efficient methods exist for computing matrix-vector product $\mathbf{G}_i \mathbf{v}, \forall i$. The latter can be easily satisfied relying on the inherit structure of $\mathbf{G}_i$.

For notational convenience, we will write $f \circ {\cal G} (\mathbf{X}) := f\big({\cal G} ( \mathbf{X}) - \mathbf{z} \big)$, where $f : \text{dom}~ f \in \mathbb{R}^d \mapsto \mathbb{R}$, and $f \circ {\cal G}: \mathbb{S}^n_+ \mapsto \mathbb{R}$. The matrix norm $\| \cdot \|$ in this section refers to Schatten $1$-norm (also known as nuclear or trace norm), and its dual norm, $\| \cdot \|_\infty$, is therefore the Schatten-$\inf$ norm (or operator norm). The inner product of matrices is standard trace inner product. We do not strict the form for vector norm.

\subsection{Memory efficient implementation of MOCO}

Applying MOCO for solving \eqref{eq.prob_mtrx} requires the storage of $n \times n$ matrices $\mathbf{X}_k$ and $\mathbf{g}_k$. Note that here $\mathbf{g}_k$ is a matrix, and we keep the same notation as Alg. \ref{alg.mcd} for consistence. This memory consumption is a significant barrier for scaling problems up. Moreover, it is extremely not economical for raised-up SDPs as discussed in Section \ref{sec.intro}. To facilitate memory efficiency in MOCO, the changes are represented below.

\textbf{Vectorized representation of $\mathbf{X}_k$.} Let $\mathbf{y}_k = {\cal G}(\mathbf{X}_k) - \mathbf{z}$. The vector $\mathbf{y}_k$ is helpful for memory saving of MOCO iterates. In particular, $\eta_k$ can be obtained using only vectors $\mathbf{y}_k$ and $\mathbf{z}$
\begin{align}
	\eta_k & = \argmin_{\eta \geq 0} f \circ {\cal G}( \eta \mathbf{X}_k) = \argmin_{\eta \geq 0} f \big( {\cal G}( \eta \mathbf{X}_k)  - \mathbf{z} \big) \nonumber \\
	& =  \argmin_{\eta \geq 0} f \big(\eta  {\cal G}( \mathbf{X}_k)  - \mathbf{z} \big) =  \argmin_{\eta \geq 0} f(\eta \mathbf{y}_k + \eta \mathbf{z} - \mathbf{z}). \nonumber  
\end{align}
Similarly, $\mathbf{y}_k$ avoids explicit use of $\mathbf{X}_k$ when solving for $\theta_k$ in Line 6 in Alg. \ref{alg.mcd_mem}. Another merit of $\mathbf{y}_k$ lies in the fact that $\nabla f \circ {\cal G}(\mathbf{X}_k) = {\cal G}^* \big(\nabla f (\mathbf{y}_k) \big)$. Owing to the linearity of ${\cal G}^*$, i.e., ${\cal G}^*(\mathbf{a} + \mathbf{b}) = {\cal G}^*(\mathbf{a}) +{\cal G}^*( \mathbf{b})$, it is possible to leverage a vector $\tilde{\mathbf{g}}_k \in \mathbb{R}^d$ to retrieve the full gradient $\mathbf{g}_k$ as ${\cal G}^*(\tilde{\mathbf{g}}_k)$; see line 4 of Alg. \ref{alg.mcd_mem}. 

\textbf{MOCO subproblem.}
The MOCO subproblem in line 5 under Schatten 1-norm is equivalent to find the minimum eigenvalue and its normalized eigenvector of ${\cal G}^*(\tilde{\mathbf{g}}_k)$. This can be carried out efficiently through shifted power method or the Lanczos method \cite{saad2011numerical}.

\textbf{Sketched representation  of $\mathbf{X}_k$.} Although $\mathbf{y}_k$ removes the explicit need of $\mathbf{X}_k$, it does not support to reconstruct $\mathbf{X}_k$. Random sketches $\mathbf{S}_k$ are adopted to address this problem in memory efficient form following \cite{yurtsever2021}. 
In particular, a random Gaussian matrix $\bm{\Omega} \in \mathbb{R}^{n \times R}$ is fixed for some predefined parameter $R \ll n$. The sketch is then defined as $\mathbf{S}_k = \mathbf{X}_k \bm{\Omega}$. The linearity of sketch also permits a simple update for $\mathbf{S}_k$ in line 8 of Alg. \ref{alg.mcd_mem}. For the ease of understanding and analyses, the update of $\mathbf{X}_k$ is written in line 9, however, this line should be omitted when coding. The overall memory consumption for Alg. \ref{alg.mcd_mem} is ${\cal O}(d + nR)$, which can be much less than ${\cal O}(n^2)$ in the naive implementation of MOCO.

\begin{algorithm}[t]
    \caption{Memory efficient MOCO for \eqref{eq.prob_mtrx}}\label{alg.mcd_mem}
    \begin{algorithmic}[1]
    	\State \textbf{Initialize:} $\mathbf{y}_0 = - \mathbf{z}$, $\delta_k = \frac{2}{k+2} \forall k, \mathbf{S}_0 = \mathbf{0} \in \mathbb{R}^{n \times r}$
			\For {$k=0,1,\dots,K$}
				\State $\eta_k = \argmin_{\eta \geq 0} f(\eta \mathbf{y}_k + \eta \mathbf{z} - \mathbf{z})$ 
				\State $\tilde{\mathbf{g}}_k = (1 - \delta_k) \tilde{\mathbf{g}}_{k-1} + \delta_k \nabla f( \eta \mathbf{y}_k + \eta \mathbf{z} - \mathbf{z} )$
				\State find $\lambda_k =\lambda_{\min} ( {\cal G}^*[ \tilde{\mathbf{g}}_k ] )$ and associated normalized eigenvector $\mathbf{q}_k$
				\State $\theta_k = \argmin_{\theta \geq 0} f(\eta_k \mathbf{y}_k + \eta_k \mathbf{z} - \mathbf{z}+\theta {\cal G}(\mathbf{q}_k\mathbf{q}_k^\top ) )$ 
				\State $\mathbf{y}_{k+1} = \eta_k \mathbf{y}_k + \eta_k \mathbf{z} - \mathbf{z}+\theta_k {\cal G}(\mathbf{q}_k\mathbf{q}_k^\top)$
				\State \textbf{Option I:} $\mathbf{S}_{k+1} = \eta_k \mathbf{S}_k + \theta_k \mathbf{q}_k (\mathbf{q}_k^\top \bm{\Omega})$
				\State \textbf{Option II:} $\mathbf{X}_{k+1} = \eta_k \mathbf{X}_k + \theta_k \mathbf{q}_k \mathbf{q}_k^\top $
				\State (optional) greedy step in Alg. \ref{alg.greedy}
			\EndFor
		\State \textbf{Return:} $\eta_{K} \mathbf{S}_{K}$ (to recover $\mathbf{X}_K$)
	\end{algorithmic}
\end{algorithm}

\textbf{Recover $\mathbf{X}_k$ from $\mathbf{S}_k$.} One can find a rank $r$ approximation $\hat{\mathbf{X}}_k$ to the real variable $\mathbf{X}_k$ using a stable implementation \cite[Algorithm 5.1]{yurtsever2021}. The reconstruction error is bounded as the following if $r < R-1$
\begin{align*}
	\mathbb{E} \big[ \| \mathbf{X}_a - \hat{\mathbf{X}}_a \| \big]	 \leq \bigg( 1 + \frac{r}{R- r - 1} \bigg) \sum_{i=r+1}^{n} \sigma_i (\mathbf{X}_a).
\end{align*}
The reconstruction error is sufficient small when the true $\mathbf{X}_k$ is low rank. This means that the memory efficiency is almost obtained for free for problems such as raised-up SDPs.

\subsection{Convergence}

 In this subsection, we will strengthen Theorems \ref{thm.primal} -- \ref{thm.stopcr} by relaxing Assumption \ref{as.2}.
The modified assumptions are listed below, and indexed with a prime to connect with its original counterpart in Section \ref{sec.assumptions}.

\begin{assumptionp}{\ref{as.1}$'$}\label{as.1p}
	$f \circ {\cal G}(\mathbf{X})$ has $L_{f \circ {\cal G}}$-Lipchitz continuous gradients, and $f(\mathbf{y})$ has $L_f$-Lipchitz continuous gradients.
\end{assumptionp}

\begin{assumptionp}{\ref{as.2}$'$}\label{as.2p}
	$f$ is strictly convex.	
\end{assumptionp}

Assumption \ref{as.2p} relaxes Assumption \ref{as.2} since it does not require strict convexity on $f \circ {\cal G}$.

Since $\mathbb{S}_n^+$ is convex, we do not state Assumption \ref{as.3} explicitly here. ${\cal K}$ and ${\cal K}^*$ are also not explicitly distinguished since $\mathbb{S}_n^+$ is self-dual. In addition to previous assumptions, it is convenient to have another bounded assumption on ${\cal G}$ or its adjoint ${\cal G}^*$. 

\begin{assumption}\label{as.4}
	The adjoint ${\cal G}^*(\cdot)$ is $\bar{G}$-Lipschitz, that is,
	$\| {\cal G}^*(\mathbf{a})\|_* \leq \bar{G} \| \mathbf{a} \|_*$. 
\end{assumption}
This is a very mild assumption in practice, and it is satisfied long as matrices $\{ \mathbf{G}_i \}_{i=1}^d$ are bounded. To see this, we have
\begin{align*}
	\| {\cal G}^*(\mathbf{a}) \|_* & \leq | a_1 | \| \mathbf{G}_1 \|_* + \ldots + | a_d| \| \mathbf{G}_d \|_* \\
	& = \langle \mathbf{m}, \tilde{\mathbf{a}} \rangle \leq \| \mathbf{m} \|  \| \mathbf{a} \|_*
\end{align*}
where $\mathbf{m}:= [\| \mathbf{G}_1 \|_*, \ldots, \| \mathbf{G}_d \|_* ]^\top$, and $\tilde{\mathbf{a}}:= [|a_1|, \ldots, |a_d| ]^\top$. The inequality above implies that $\bar{G} \leq \|\mathbf{m} \|$.

\begin{theoremp}{\ref{thm.primal}$'$}\label{thm.primal'}
(Primal convergence.)
Under Assumptions \ref{as.1p} and \ref{as.2p}, Alg. \ref{alg.mcd_mem} for solving \eqref{eq.prob_mtrx} ensures that
\begin{align*}
	f(\eta_k\mathbf{y}_k ) = f \circ {\cal G}(\eta_k\mathbf{X}_k)   \leq f \circ {\cal G}(\mathbf{X}^*)  + \frac{2L_{f \circ {\cal G}} \|\mathbf{X} \|^2}{k+1}.	
\end{align*}
Moreover, let $\Psi^*$ collect the optimal solutions of \eqref{eq.prob_mtrx}, and $\hat{\mathbf{X}}_k$ denote the reconstructed matrix from $\eta_k\mathbf{S}_k$, then the reconstruction error satisfies that 
\begin{align*}
	\limsup_{k  \rightarrow \infty} \mathbb{E}[\text{dist}(\hat{\mathbf{X}}_k, \Psi^*) ] \leq \bigg( 1 + \frac{r}{R- r - 1} \bigg) \max_{\mathbf{X}^* \in \Psi^*} \Sigma_r(\mathbf{X}^*)
\end{align*}
where $\Sigma_r(\mathbf{X}^*):=\sum_{i=r+1}^{n} \sigma_i (\mathbf{X}^*)$.
\end{theoremp}

Theorem \ref{thm.primal'} establishes that the convergence of $\mathbf{y}_k$, together with a bound for on the reconstruction error. In particular, if the desirable $\mathbf{X}^*$ is indeed low rank (smaller than $r$), the reconstruction error is $0$. This once again suggests that the memory efficient MOCO is suitable for problems with low rank solutions.

\begin{theoremp}{\ref{thm.dual}$'$}\label{thm.dual'}
(Dual convergence.)
Under Assumptions \ref{as.1p}, \ref{as.2p}, and \ref{as.4}, Alg. \ref{alg.mcd_mem} for solving \eqref{eq.prob_mtrx} ensures that
	\begin{align*}
	\big[ \text{dist}^*(\nabla f \circ {\cal G}(\eta_k\mathbf{X}_k), \mathbb{S}_n^+) \big]^2 \leq \frac{ {4 \bar{G}^2 L_f L_{f \circ {\cal G}} \| \mathbf{X}^* \|^2 } }{k+1}.
	\end{align*}
	Furthermore, let $\hat{\mathbf{X}}_k$ denote the reconstructed matrix from $\eta_k\mathbf{S}_k$, then the reconstruction error satisfies that 
	\begin{align*}
		&\limsup_{k  \rightarrow \infty} \mathbb{E}[\text{dist}^*(\nabla f \circ {\cal G}(\hat{\mathbf{X}}_k) , \mathbb{S}_n^+) ] \\
		& ~~~~~~~ \leq \bigg( 1 + \frac{r}{R- r - 1} \bigg) L_{f\circ {\cal G}}\max_{\mathbf{X}^* \in \Psi^*} \Sigma_r(\mathbf{X}^*)
	\end{align*}
	where $\Sigma_r(\mathbf{X}^*)$ is defined the same as Theorem \ref{thm.primal'}. 
\end{theoremp}



\begin{theoremp}{\ref{thm.stopcr}$'$}\label{thm.stopcr'}
(Stopping criterion.)
Under Assumptions \ref{as.1p}, \ref{as.2p}, and \ref{as.4}, Alg. \ref{alg.mcd_mem} for solving \eqref{eq.prob_mtrx} ensures that
\begin{align*}
	\big[\text{dist}_*({\cal G}^*(\tilde{\mathbf{g}}_k), \mathbb{S}_n^+)\big]^2 \leq {\cal O}\bigg(\frac{\bar{G}^2 L_f L_{f \circ {\cal G}} \|  \mathbf{X}^*\|^2}{k} \bigg)
\end{align*}
\end{theoremp}

Theorem \ref{thm.stopcr'} suggests a natural stopping criterion to certify that $\eta_k\mathbf{X}_k$ is near optimal, which is
\begin{align}
	\lambda_k \geq - \sqrt{\epsilon}.
\end{align}
In other words, MOCO can be terminated long as ${\cal G}^*( \tilde{\mathbf{g}}_k ) $ is almost positive semidefinite. Note that the stopping criterion is irrelevant to the reconstruction error, because it depends on neither $\mathbf{S}_k$ nor $\mathbf{X}_k$.

Going beyond theories, two heuristics are introduced in the following subsections to further improve the practical merits of MOCO.

\subsection{Practical heuristic 1: Greedy step}

\begin{algorithm}[t]
    \caption{Greedy step at iteration $k$}\label{alg.greedy}
    \begin{algorithmic}[1]
    	\State find $(t_k, \mathbf{U}_k)$ by solving \eqref{eq.greedy} 
    	\State $\mathbf{y}_{k+1} \leftarrow  t_k^2(\mathbf{y}_{k+1} + \mathbf{z}) + {\cal G} ( \mathbf{U}_k\mathbf{U}_k^\top) - \mathbf{z} $
    	\State $\mathbf{S}_{k+1} \leftarrow  t_k^2(\mathbf{S}_{k+1} + \mathbf{z}) + \mathbf{U}_k (\mathbf{U}_k^\top \bm{\Omega}) - \mathbf{z} $
    \end{algorithmic}
\end{algorithm}

\begin{figure*}[t]
	\centering
	\begin{tabular}{c}
		\includegraphics[width=0.97\textwidth]{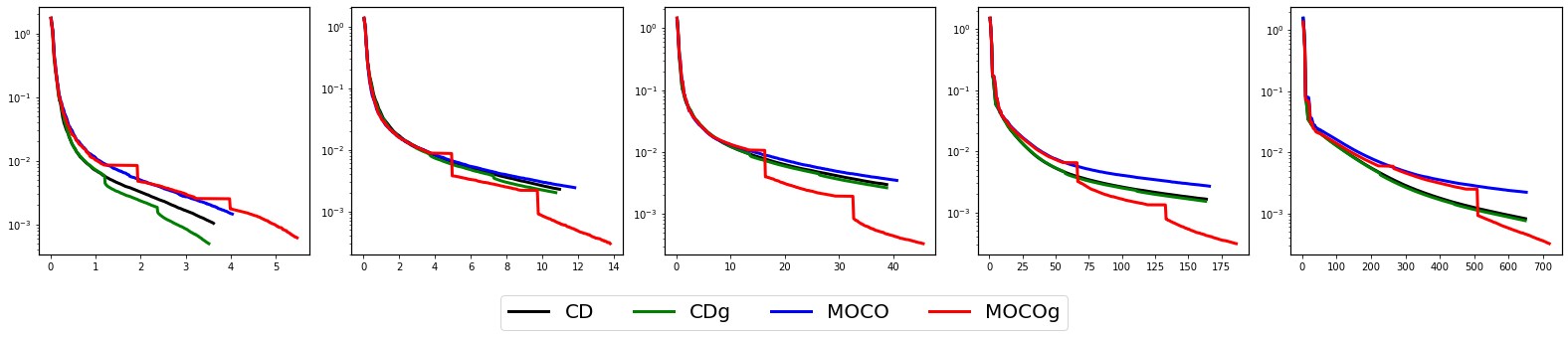} \\
	\end{tabular}
	\caption{Performances (runtime vs primal error) of different algorithms for the matrix completion problem \eqref{eq.mtrx_complete_relax}. From left to right, the sizes of problems are $n= 100, 200, 400, 800, 1600$.}
	 \label{fig.mtrx_comp}
\end{figure*}

This heuristic aims to handle raised-up SDPs through Burer-Monteiro (BM) factorization \cite{burer2003, bhojanapalli2016}. The idea is to greadily move $\mathbf{X}_{k+1}$ to a point on another ray to reduce the objective value. This can be done by solving the following unconstrained problem through any descent method 
\cite{shewchuk1994introduction}
\begin{align}\label{eq.greedy}
	(t_k, \mathbf{U}_k) & = \argmin_{t \in \mathbb{R}, \mathbf{U} \in \mathbb{R}^{n \times r}} f\big({\cal G} ( t^2 \mathbf{X}_{k+1} + \mathbf{U}\mathbf{U}^\top) - \mathbf{z} \big) 
	\\
	& = \argmin_{t , \mathbf{U} } f\big( t^2(\mathbf{y}_{k+1} + \mathbf{z}) + {\cal G} ( \mathbf{U}\mathbf{U}^\top) - \mathbf{z} \big).  \nonumber
\end{align}
Note that $t_k \mathbf{X}_k + \mathbf{U}_k \mathbf{U}_k^\top$ is a positive semidefinite matrix. Then, the feasible point $t_k \mathbf{X}_k + \mathbf{U}_k \mathbf{U}_k^\top$ is used as the starting point of next iteration of Alg. \ref{alg.mcd_mem}. Accordingly, the way to update $\mathbf{y}_k$ and $\mathbf{S}_k$ based on \eqref{eq.greedy} is given in Alg. \ref{alg.greedy}. The greedy step is optional and can be helpful to run every a few (e.g., $100$) iterations to speedup convergence.

Another manner to understand the greedy step is by viewing MOCO as a theoretical justified wrapper for the BM approach, where the convergence of latter is difficult to establish. Because problem \eqref{eq.greedy} is solved through a \textit{descent} approach, the greedy step aided MOCO converges naturally. Note that when choosing a proper solver for \eqref{eq.greedy}, the memory consumption of the greedy step aided MOCO is still ${\cal O}(d + nR)$. Although the greedy step also applies to CD \cite{duchi2020conic}, we find that it is more suitable for MOCO because of the improved performance as shown later in our numerical tests.

\subsection{Practical heuristic 2: magical $\theta_k$} 

Next, we introduce another practical variant when $\| \mathbf{X}^*\|$ can be estimated. This variant can be useful for raised-up SDPs especially when $ \mathbf{X}^* = \mathbf{x}^* (\mathbf{x}^*)^\top $ for some vector $\mathbf{x}^*$. In this case, it can be possible to use the relation $\| \mathbf{X}^*\| = \| \mathbf{x}^* \|_2^2$ to estimate $\| \mathbf{X}^*\|$. An example will be provided in Section \ref{sec.num.phase} together with numerical tests.

This heuristic is motivated by the empirical wisdom that line search can be conservative for numerical performances of heavy ball momentum for FW \cite{li2021heavy}. Let $M > 0$ be an estimate of $\| \mathbf{X}^*\|$, then our heuristic step size is $\theta_k = \frac{2 M}{k + 2}$. This step size comes from the detailed derivation of Theorem \ref{thm.primal}; see the first line of Appendix \ref{apdx.thm.primal}. For problems where $M$ is difficult to estimate, it is also possible we can run MOCO for a few iterations, then use $\| \mathbf{X}_k \| $ as an estimate of $\| \mathbf{X}^* \| $. The heuristic $\theta_k$ eliminates the need for line search, therefore saving runtime.

\section{Numerical tests}

Experiments with synthetic and real data are conducted to visualize the performance of the proposed MOCO and its practical heuristics.

\subsection{Matrix completion}

MOCO is first tested on matrix completion problems using synthetic data. Suppose the ground truth matrix $\mathbf{A} \in \mathbb{S}^n_+$ to be recovered is low rank and positive semidefinite. We are given noisy entries of $\mathbf{A}$ sampled randomly, that we denote as $b_{ij} = A_{ij} + \epsilon_{ij}$ for some index $(i,j) \in {\cal I}$, where $\epsilon_{ij}$ are zero mean i.i.d. Gaussian random variables. Let $\mathbf{A} = \mathbf{V} \mathbf{V}^\top$ for some $\mathbf{V} \in \mathbb{R}^{n \times 3}$ denote the low-rank ground truth. The estimated matrix can be found by solving the following problem
\begin{align}\label{eq.mtrx_complete_relax}
	\min_{\mathbf{X}}  ~~ \frac{1}{2} \sum_{(i,j) \in {\cal I}} (X_{ij} - b_{ij})^2 ~~~~
		\text{s.t.}  ~~ \mathbf{X} \in \mathbb{S}_n^+ .
\end{align}

To understand how MOCO scales, we consider \eqref{eq.mtrx_complete_relax} with number of data $n \in \{ 100, 200, 400, 800, 1600 \}$. Following \cite{duchi2020conic}, we sample every entry in the upper left $10\times 10$ block and other entries with probability $0.1$. The Gaussian noise $\epsilon_{ij}$ is randomly generated so that the SNR is $20$dB.

The benchmark algorithms are chosen as CD and CD with a greedy heuristic (CDg) \cite{duchi2020conic}. Our numerical tests rely on memory efficient implementation of MOCO. MOCO with the greedy heuristic (MOCOg) is also considered to improve numerical performance. Each tested algorithm is run for $300$ iterations. The primal error versus runtime is plotted in Figure \ref{fig.mtrx_comp}. For the matrix completion problem, MOCO exhibits slightly worse performance compared to CD. On the other hand, the greedy step appears to be more suitable for MOCO since it clearly boosts the performance of MOCO but not CD with the only exception on the test case with smallest scale $n=100$. In addition, given the same amount of time, MOCOg achieves the lowest primal error compared with other tested algorithms, thus confirming its scalability and efficiency. The greedy step does not make enough progress for CD, but it significantly helps MOCO. This empirically suggests that the merits of the greedy step are amplified by the momentum in MOCO.

\subsection{Phase retrieval}\label{sec.num.phase}

\begin{figure*}[t]
	\centering
	\begin{tabular}{c}
		\hspace{-1.cm}
		\includegraphics[width=1.07\textwidth]{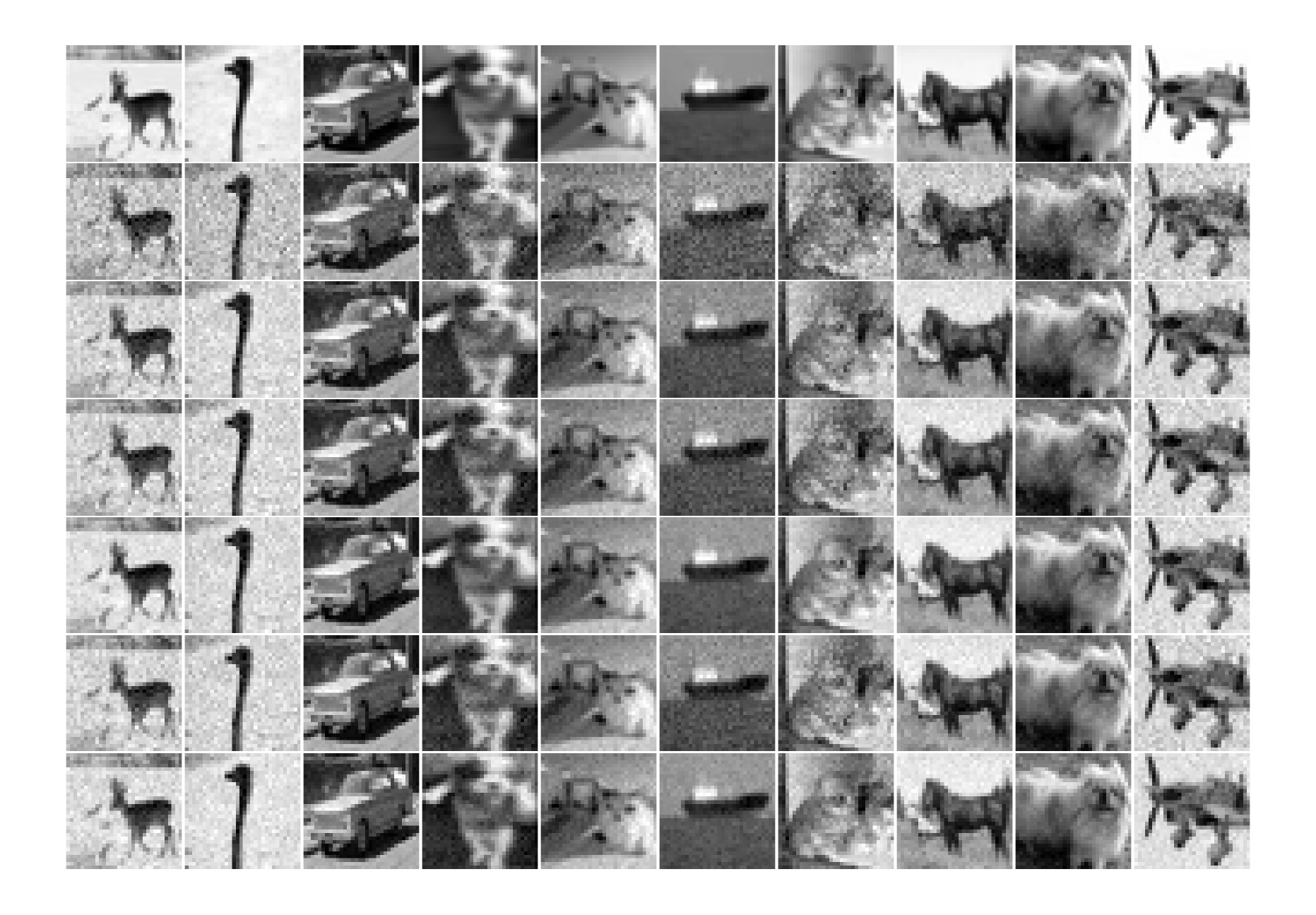} \\
		\hspace{-1.1cm}
		\includegraphics[width=0.97\textwidth]{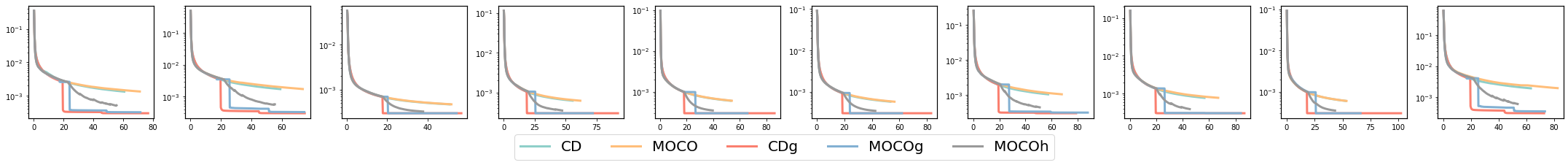}	
	\end{tabular}
	\vspace{-0.3cm}
	\caption{Performances of various algorithms for the phase retrieval problem. Each column contains the result using a specific image. The first row plots raw images, and other rows (from 2nd  to 7th) contain images recovered using FW, CD, MOCO, CDg, MOCOg, MOCOh, respectively. And the last row lists the optimality error vs iteration of compared approaches.}
	 \label{fig.phase_ret}
\end{figure*}

Suppose that $\mathbf{x} \in \mathbb{R}^n$ is a signal to be retrieved from measurements $b_i = (\mathbf{a}_i^\top \mathbf{x})^2 + \epsilon_i$, where $\mathbf{a}_i \in \mathbb{R}^n$ are rows of matrix $\mathbf{A} = [\mathbf{D}\mathbf{S}_1, \ldots, \mathbf{D}\mathbf{S}_m]^\top$ with $\mathbf{D}$ being the discrete cosine transform and $\mathbf{S}_1, \ldots, \mathbf{S}_m$ being diagonal matrices of independent random signs. One means to recover the original signal is to solve the following problem 
\begin{align*}
	\min_{\mathbf{x} \in \mathbb{R}^n} \frac{1}{mn} \sum_{i=1}^{mn} \| b_i - (\mathbf{a}_i^\top \mathbf{x})^2 \|_2^2 + \gamma \| \mathbf{x} \|_2^2.
\end{align*}
The unsatisfactory of this formulation resides in the fact that the objective function is a polynomial of forth order. This challenges optimization since it is non-smooth. Raised-up SDP is the remedy. Noticing that $(\mathbf{a}_i^\top \mathbf{x})^2 = \mathbf{a}_i^\top \mathbf{x} \mathbf{x}^\top \mathbf{a}_i := \mathbf{a}_i^\top \mathbf{X} \mathbf{a}_i$ where $\mathbf{X}:= \mathbf{x}\mathbf{x}^\top$, we can reformulate the problem as
\begin{align}\label{eq.phase_ret}
	\min_{\mathbf{X}\in \mathbb{S}_+^n } ~  \frac{1}{mn} \big\| {\cal A} ( \mathbf{X}) - \mathbf{b} \big\|^2 + \gamma \tr(\mathbf{X})
\end{align}
where ${\cal A} = [\tr(\mathbf{a}_1 \mathbf{a}_1^\top \mathbf{X}), \ldots, \tr(\mathbf{a}_{mn} \mathbf{a}_{mn}^\top \mathbf{X})]$. Since $\mathbf{X} = \mathbf{x}\mathbf{x}^\top$, it is natural to assume that there exists a rank-1 optimal solution $\mathbf{X}^* = \mathbf{x}^* (\mathbf{x}^*)^\top$. The rank-$1$ assumption also enables an estimate of $\| \mathbf{X}^* \|$ to use in the heuristic $\theta_k$
\begin{align}
	\| \mathbf{X}^* \| \!= \!\tr (\mathbf{X}^*) = \| \mathbf{x}^* \|_2^2  = \frac{1}{m} \sum_{i=1}^{mn} (\mathbf{a}_i^\top \mathbf{x}^*)^2 \!\approx\! \frac{1}{m} \! \sum_{i=1}^{mn} b_i.
\end{align}

For the experiment setup, $10$ images from CIFAR10 dataset\footnote{https://www.cs.toronto.edu/~kriz/cifar.html} are randomly chosen as the raw signal $\mathbf{x}^*$. The Gaussian noise $\epsilon_i$ is generated with $20$dB SNR. Other parameters are set to $m=10$ and $\gamma = 5 \times 10^{-5}$. The benchmark algorithms are chosen as FW, CD, and CD with greedy heuristic (CDg). When working with FW, we add another constraint $ \tr(\mathbf{X}) \leq \frac{2}{m} \sum_{i=1}^{mn} b_i$ to ensure the compactness of the constraint set, where the right hand side of this inequality is roughly $2\|\mathbf{X}^* \|$. Three MOCO variants are tested: MOCO in Alg. \ref{alg.mcd_mem}, MOCO with greedy heuristic (MOCOg), and MOCO with heuristic $\theta_k$ (MOCOh). All algorithms are run for $300$ iterations. We use $R=3$ for the sketches.

The original and recovered figures are shown in Fig. \ref{fig.phase_ret}, where the first row lists raw images, and other rows are recovered images via FW, CD, MOCO, CDg, MOCOg, and MOCOh, respectively. Among all implemented approaches, the recovered images using FW have the worst quality. CD and MOCO have almost the same recovery quality, and CDg, MOCOg and MOCOh share the best figure quality. This demonstrates that MOCO and CD not only improve numerical performances over FW, but remove the need for the compact domain requirement. 

To further showcase the merits of MOCO over CD, we also plot $f \circ {\cal G} (\mathbf{X}_k) - f \circ {\cal G}(\mathbf{X}^*)$ versus runtime. The loss curve for FW is omitted here because it works on a different problem from \eqref{eq.phase_ret} due to the added constraint. Despite the runtime of MOCO is longer than CD because of updating $\tilde{\mathbf{g}}_k$, the runtime of MOCOg is less than that of CDg. Hence, the greedy heuristic is better use with MOCO than CD. Although we do not have an exact explanation, our guess is that the loss curvature of the greedy subproblem could be ill-conditioned in CD. In addition, although relying on a heuristic $\theta_k$, MOCOh often converges faster than MOCO, and even matches to the performance of MOCOg sometimes. As the heuristic $\theta_k$ eliminates the need for line search, the runtime of MOCOh is shorter than MOCO.

\section{Conclusion}

This paper revisits conic descent (CD) for conic programming problems. CD is refined through a geometrical interpretation that has matching mathematical foundation in the dual domain. Then a new approach, \underline{mo}mentum \underline{co}nic descent (MOCO), is proposed to improve CD empirically and theoretically. Lastly, the dual behavior of MOCO (as well as CD) is comprehensively examined, where it is discusses about stopping criterion and opportunities to accelerate convergence via preconditioning. A memory efficient implementation of MOCO for SDPs is then developed based. Memory efficiency is achieved almost for free given the low rankness of the solution. Numerical results further validate the efficiency of proposed MOCO and its practical variants.

\bibliographystyle{IEEEtranS}
\bibliography{myabrv,datactr}

\begin{thebibliography}{10}
\providecommand{\url}[1]{#1}
\csname url@samestyle\endcsname
\providecommand{\newblock}{\relax}
\providecommand{\bibinfo}[2]{#2}
\providecommand{\BIBentrySTDinterwordspacing}{\spaceskip=0pt\relax}
\providecommand{\BIBentryALTinterwordstretchfactor}{4}
\providecommand{\BIBentryALTinterwordspacing}{\spaceskip=\fontdimen2\font plus
\BIBentryALTinterwordstretchfactor\fontdimen3\font minus
  \fontdimen4\font\relax}
\providecommand{\BIBforeignlanguage}[2]{{%
\expandafter\ifx\csname l@#1\endcsname\relax
\typeout{** WARNING: IEEEtranS.bst: No hyphenation pattern has been}%
\typeout{** loaded for the language `#1'. Using the pattern for}%
\typeout{** the default language instead.}%
\else
\language=\csname l@#1\endcsname
\fi
#2}}
\providecommand{\BIBdecl}{\relax}
\BIBdecl

\bibitem{mosek}
``{MOSEK} modeling cookbook 3.3.0,''
  \url{https://docs.mosek.com/modeling-cookbook/powo.html}.

\bibitem{bhojanapalli2016}
S.~Bhojanapalli, A.~Kyrillidis, and S.~Sanghavi, ``Dropping convexity for
  faster semi-definite optimization,'' in \emph{Conference on Learning
  Theory}.\hskip 1em plus 0.5em minus 0.4em\relax PMLR, 2016, pp. 530--582.

\bibitem{boyd2004}
S.~Boyd, S.~P. Boyd, and L.~Vandenberghe, \emph{Convex optimization}.\hskip 1em
  plus 0.5em minus 0.4em\relax Cambridge university press, 2004.

\bibitem{burer2003}
S.~Burer and R.~D. Monteiro, ``A nonlinear programming algorithm for solving
  semidefinite programs via low-rank factorization,'' \emph{Mathematical
  Programming}, vol.~95, no.~2, pp. 329--357, 2003.

\bibitem{cui2020projecting}
Y.~Cui, L.~Liang, D.~Sun, and K.-C. Toh, ``Projecting onto the degenerate
  doubly nonnegative cone,'' \emph{arXiv preprint arXiv:2009.11272}, 2020.

\bibitem{duchi2020conic}
J.~C. Duchi, O.~Hinder, A.~Naber, and Y.~Ye, ``Conic descent and its
  application to memory-efficient optimization over positive semidefinite
  matrices,'' \emph{Advances in Neural Information Processing Systems},
  vol.~33, pp. 8308--8317, 2020.

\bibitem{dur2010}
M.~D{\"u}r, ``Copositive programming--a survey,'' in \emph{Recent advances in
  optimization and its applications in engineering}.\hskip 1em plus 0.5em minus
  0.4em\relax Springer, 2010, pp. 3--20.

\bibitem{frank1956}
M.~Frank and P.~Wolfe, ``An algorithm for quadratic programming,'' \emph{Naval
  Research Logistics Quarterly}, vol.~3, no. 1-2, pp. 95--110, 1956.

\bibitem{hajek2016}
B.~Hajek, Y.~Wu, and J.~Xu, ``Achieving exact cluster recovery threshold via
  semidefinite programming,'' \emph{IEEE Transactions on Information Theory},
  vol.~62, no.~5, pp. 2788--2797, 2016.

\bibitem{hanasusanto2018conic}
G.~A. Hanasusanto and D.~Kuhn, ``Conic programming reformulations of two-stage
  distributionally robust linear programs over wasserstein balls,''
  \emph{Operations Research}, vol.~66, no.~3, pp. 849--869, 2018.

\bibitem{harchaoui2015}
Z.~Harchaoui, A.~Juditsky, and A.~Nemirovski, ``Conditional gradient algorithms
  for norm-regularized smooth convex optimization,'' \emph{Mathematical
  Programming}, vol. 152, no.~1, pp. 75--112, 2015.

\bibitem{jaggi2013}
M.~Jaggi, ``Revisiting {F}rank-{W}olfe: {P}rojection-free sparse convex
  optimization.'' in \emph{Proc. Intl. Conf. on Machine Learning}, 2013, pp.
  427--435.

\bibitem{kakade2009duality}
S.~Kakade, S.~Shalev-Shwartz, A.~Tewari \emph{et~al.}, ``On the duality of
  strong convexity and strong smoothness: Learning applications and matrix
  regularization.''

\bibitem{kato2007}
T.~Kato, H.~Kashima, M.~Sugiyama, and K.~Asai, ``Multi-task learning via conic
  programming,'' \emph{Advances in Neural Information Processing Systems},
  vol.~20, 2007.

\bibitem{li2020}
B.~Li, M.~Coutino, G.~B. Giannakis, and G.~Leus, ``A momentum-guided
  {F}rank-{W}olfe algorithm,'' \emph{IEEE Trans. on Signal Processing},
  vol.~69, pp. 3597--3611, 2021.

\bibitem{li2021heavy}
B.~Li, A.~Sadeghi, and G.~Giannakis, ``Heavy ball momentum for conditional
  gradient,'' \emph{Proc. Advances in Neural Info. Process. Syst.}, vol.~34,
  2021.

\bibitem{li2020extra}
B.~Li, L.~Wang, G.~B. Giannakis, and Z.~Zhao, ``Enhancing {F}rank {W}olfe with
  an extra subproblem,'' in \emph{Proc. of AAAI Conf. on Artificial
  Intelligence}, 2021.

\bibitem{locatello2017greedy}
F.~Locatello, M.~Tschannen, G.~R{\"a}tsch, and M.~Jaggi, ``Greedy algorithms
  for cone constrained optimization with convergence guarantees,''
  \emph{Advances in Neural Information Processing Systems}, vol.~30, 2017.

\bibitem{nesterov2000semidefinite}
Y.~Nesterov, H.~Wolkowicz, and Y.~Ye, ``Semidefinite programming relaxations of
  nonconvex quadratic optimization,'' in \emph{Handbook of semidefinite
  programming}.\hskip 1em plus 0.5em minus 0.4em\relax Springer, 2000, pp.
  361--419.

\bibitem{nesterov2004}
Y.~Nesterov, \emph{Introductory lectures on convex optimization: A basic
  course}.\hskip 1em plus 0.5em minus 0.4em\relax Springer Science \& Business
  Media, 2004, vol.~87.

\bibitem{saad2011numerical}
Y.~Saad, \emph{Numerical methods for large eigenvalue problems: revised
  edition}.\hskip 1em plus 0.5em minus 0.4em\relax SIAM, 2011.

\bibitem{shewchuk1994introduction}
J.~R. Shewchuk \emph{et~al.}, ``An introduction to the conjugate gradient
  method without the agonizing pain,'' 1994.

\bibitem{vandenberghe1996}
L.~Vandenberghe and S.~Boyd, ``Semidefinite programming,'' \emph{SIAM review},
  vol.~38, no.~1, pp. 49--95, 1996.

\bibitem{wang2015robust}
X.~Wang, Y.~Zhang, G.~B. Giannakis, and S.~Hu, ``Robust smart-grid-powered
  cooperative multipoint systems,'' \emph{IEEE Transactions on Wireless
  Communications}, vol.~14, no.~11, pp. 6188--6199, 2015.

\bibitem{yu2004iterative}
W.~Yu, W.~Rhee, S.~Boyd, and J.~M. Cioffi, ``Iterative water-filling for
  gaussian vector multiple-access channels,'' \emph{IEEE Transactions on
  Information Theory}, vol.~50, no.~1, pp. 145--152, 2004.

\bibitem{yurtsever2021}
A.~Yurtsever, J.~A. Tropp, O.~Fercoq, M.~Udell, and V.~Cevher, ``Scalable
  semidefinite programming,'' \emph{SIAM Journal on Mathematics of Data
  Science}, vol.~3, no.~1, pp. 171--200, 2021.

\bibitem{zhang2021}
Y.~Zhang, B.~Li, and G.~B. Giannakis, ``Accelerating frank-wolfe with weighted
  average gradients,'' in \emph{ICASSP 2021-2021 IEEE International Conference
  on Acoustics, Speech and Signal Processing (ICASSP)}.\hskip 1em plus 0.5em
  minus 0.4em\relax IEEE, 2021, pp. 5529--5533.

\end{thebibliography}


\appendix

\subsection{Proof of Lemma \ref{lemma.raymin}}
\begin{proof}
	Following line 3 of MOCO, this lemma can be verified by considering two cases: i) $\eta_k = 0$, and, ii) $\eta_k \neq 0$. For case i) equation \eqref{eq.zero} holds directly, while for ii), we have from optimality condition that $\langle \mathbf{x}_k, \nabla f(\eta_k \mathbf{x}_k) \rangle = 0$. 
\end{proof}

\subsection{Proof of Lemma \ref{lemma.phi}}\label{apdx.lemma.phi}

To prove i), let ${\cal X}_r:= \{\mathbf{x}| \mathbf{x} \in {\cal K}, \| \mathbf{x} \| \leq r \}$. By rewriting $\Phi_{k+1}(\mathbf{x})$, it is not hard to see that
		\begin{align*}
			& \min_{\mathbf{x} \in {\cal X}_r} \Phi_{k+1}(\mathbf{x}) \\
			 \Leftrightarrow  & \min_{\mathbf{x} \in {\cal X}_r} (1 - \delta_k) \Phi_k(\mathbf{x}) + \delta_k  \langle \nabla f(\eta_k\mathbf{x}_k), \mathbf{x} \rangle  \\
			 \Leftrightarrow & \min_{\mathbf{x} \in {\cal X}_r} \langle \mathbf{g}_k, \mathbf{x} \rangle.
		\end{align*}
		Given line 5 of Alg. \ref{alg.mcd}, we have that $\mathbf{v}_k$ minimizes $\Phi_{k+1}(\mathbf{x})$ over ${\cal X}_1$. This completes the proof of property i).
		
		To see ii), we have that
		 \begin{align}\label{eq.phi_f}
			&~~~~~ \Phi_{k+1}( \mathbf{x}) \\
			 & = (1-\delta_k) \Phi_k(\mathbf{x} ) + \delta_k \big[ f(\eta_k\mathbf{x}_k) + \langle \nabla f(\eta_k\mathbf{x}_k), \mathbf{x}  \rangle \big] \nonumber \\
			& \stackrel{(a)}{=} (1-\delta_k) \Phi_k(\mathbf{x}) + \delta_k \big[ f(\eta_k\mathbf{x}_k) + \langle \nabla f(\eta_k\mathbf{x}_k), \mathbf{x} - \eta_k \mathbf{x}_k \rangle \big] \nonumber \\
			& \stackrel{(b)}{\leq} (1-\delta_k) \Phi_k(\mathbf{x}) +  \delta_k \big( f(\mathbf{x}) - \xi_k (\mathbf{x}) \big) \nonumber \\
			& \stackrel{(c)}{\leq} f(\mathbf{x}) - \rho_k (\mathbf{x}) \nonumber 
		\end{align}
		where (a) is because of \eqref{eq.zero}; (b) is by strict convexity of $f(\mathbf{x})$, and $\xi_k(\mathbf{x}) \geq 0$ with equation holding only at $\mathbf{x} = \eta_k\mathbf{x}_k$. Inequality (c) is obtained by unrolling $ (1-\delta_k) \Phi_k(\mathbf{x}) $, and $\rho_k(\mathbf{x})\geq 0$ is weighted average of $\{ \xi_\tau(\mathbf{x}) \}_{\tau=0}^k$.

		 For notational convenience, let ${\cal X}_*:= {\cal X}_{\| \mathbf{x}^* \|}$. We then have that $\| \mathbf{x}^* \| \mathbf{v}_k$ minimizes $\Phi_{k+1}(\mathbf{x})$ over ${\cal X}_*$. This is because that $\min_{\mathbf{x} \in {\cal X}_*} \Phi_{k+1}(\mathbf{x}) \Leftrightarrow  \min_{\mathbf{x} \in {\cal X}_*} \langle \mathbf{g}_k, \mathbf{x} \rangle $, where the second optimization problem has a minimal smaller than $0$ since we must have ${\mathbf{0}} \in {\cal X}_*$. Consider two cases: i) $\min_{\mathbf{x} \in {\cal X}_1} \langle \mathbf{g}_k, \mathbf{x} \rangle $ (i.e., line 5) is minimized by $\mathbf{v}_k = \mathbf{0}$, in this case $\| \mathbf{x}^* \| \mathbf{v}_k = \mathbf{0}$; and ii) $\min_{\mathbf{x} \in {\cal X}_1} \langle \mathbf{g}_k, \mathbf{x} \rangle  = \langle \mathbf{g}_k, \mathbf{v}_k \rangle < 0 $, where the smallest value of $\min_{\mathbf{x} \in {\cal X}_*} \langle \mathbf{g}_k, \mathbf{x} \rangle $ must be achieved by $\| \mathbf{x}^* \| \mathbf{v}_k$ since one can rewrite  $\min_{\mathbf{x} \in {\cal X}_*} \langle \mathbf{g}_k, \mathbf{x} \rangle $ as
		\begin{align*}
			\min_{\mathbf{x}}~~ & \| \mathbf{x}_* \|  \Big\langle  \mathbf{g}_k, \frac{\mathbf{x}}{\| \mathbf{x}^*\|}  \Big\rangle \\
			\text{s.t.} &~~ \frac{\mathbf{x}}{\| \mathbf{x}^*\|} \in {\cal K} \\
			& ~~ \Big\| \frac{\mathbf{x}}{\| \mathbf{x}^*\|} \Big\|  \leq 1
		\end{align*}
		where the first constraint is because of the definition of cone. Letting $\mathbf{y}:=  \frac{\mathbf{x}}{\| \mathbf{x}^*\|}$, it can be seen that  $\| \mathbf{x}^* \| \mathbf{v}_k$ indeed is a minimizer for this problem.
		 
		Since $\| \mathbf{x}^* \| \mathbf{v}_k$ minimizes $\Phi_{k+1}(\mathbf{x})$ over ${\cal X}_*$ and $\mathbf{x}^* \in {\cal X}_*$, we have
		\begin{align}
			\Phi_{k+1}(\|\mathbf{x}^* \|\mathbf{v}_k) + \rho_k(\mathbf{x}^*)  & \leq \Phi_{k+1}(\mathbf{x}^*) + \rho_k(\mathbf{x}^*) \leq f(\mathbf{x}^*)
		\end{align}
		where the last inequality is because of \eqref{eq.phi_f}. Let $\rho_k:= \rho_k(\mathbf{x}^*) $, and we clearly have $\rho_k = 0$ only if $\{\eta_\tau \mathbf{x}_\tau \equiv \mathbf{x}^* \}_{\tau = 0}^k$. The lemma is thus proved.

\subsection{Proof of Theorem \ref{thm.primal}}\label{apdx.thm.primal}
Let $\tilde{\mathbf{x}}_{k+1} = \eta_k\mathbf{x}_k + \tilde{\theta}_k \mathbf{v}_k $ with $\tilde{\theta}_k =  \delta_k \| \mathbf{x}^* \|$. We must have $f(\mathbf{x}_{k+1})  \leq f(\tilde{\mathbf{x}}_{k+1})$ since $\theta_k$ is obtained via line search as in line 6. Given smoothness, we then have
		\begin{align}\label{eq.smooth}
			&~~~~~ f(\mathbf{x}_{k+1})  \leq f(\tilde{\mathbf{x}}_{k+1})\\
			& \leq f(\eta_k\mathbf{x}_k) + \langle \nabla f(\eta_k\mathbf{x}_k) , \tilde{\mathbf{x}}_{k+1} - \eta_k\mathbf{x}_k \rangle + \frac{L}{2} \| \tilde{\mathbf{x}}_{k+1} - \eta_k\mathbf{x}_k \|^2 \nonumber  \\
			& = f(\eta_k\mathbf{x}_k) + \tilde{\theta}_k \langle \nabla f(\eta_k\mathbf{x}_k) ,  \mathbf{v}_k \rangle + \frac{L}{2} \| \tilde{\theta}_k\mathbf{v}_k \|^2 \nonumber \\
			& \leq f(\eta_k\mathbf{x}_k) + \tilde{\theta}_k \langle \nabla f(\eta_k\mathbf{x}_k) ,  \mathbf{v}_k \rangle + \frac{L \tilde{\theta}_k^2 }{2} \nonumber 
		\end{align}	
		where the last inequality is because $\|  \mathbf{v}_k\|\leq 1$ as in line 5 of Alg. \ref{alg.mcd}. On the other hand, we have
		\begin{align}\label{eq.phi_k_relation}
			& ~~~~ \Phi_{k+1}( \| \mathbf{x}^* \| \mathbf{v}_k) \\
			& = (1-\delta_k) \Phi_k(\| \mathbf{x}^* \| \mathbf{v}_k)  \nonumber \\
			& ~~~~~~~~~~~~~~~~~ + \delta_k \big[ f(\eta_k\mathbf{x}_k) + \langle \nabla f(\eta_k\mathbf{x}_k), \| \mathbf{x}^* \| \mathbf{v}_k  \rangle \big] \nonumber  \\
			& \stackrel{(a)}{\geq} (1-\delta_k) \Phi_k(\| \mathbf{x}^* \| \mathbf{v}_{k-1})  \nonumber \\
			&  ~~~~~~~~~~~~~~~~~~~~~ + \delta_k \big[ f(\eta_k\mathbf{x}_k) + \langle \nabla f(\eta_k\mathbf{x}_k), \| \mathbf{x}^* \| \mathbf{v}_k  \rangle \big] \nonumber
		\end{align}
		where (a) is because that $\| \mathbf{x}^* \| \mathbf{v}_{k-1}$ minimizes $\Phi_k(\mathbf{x})$ over $\{\mathbf{x}| \mathbf{x} \in {\cal K}, \| \mathbf{x} \| \leq \| \mathbf{x} \|^* \}$. Combining \eqref{eq.phi_k_relation} and \eqref{eq.smooth}, we have
		\begin{align}
			& ~~~~f(\mathbf{x}_{k+1}) - \Phi_{k+1}( \| \mathbf{x}^* \| \mathbf{v}_k) \\
			& \leq   (1-\delta_k)\big[ f(\eta_k\mathbf{x}_k) - \Phi_k(\| \mathbf{x}^* \| \mathbf{v}_{k-1}) \big]  \nonumber \\
			& ~~~~~~~ + \big(\tilde{\theta}_k - \delta_k \| \mathbf{x}^* \|  \big) \langle \nabla f(\eta_k\mathbf{x}_k) ,  \mathbf{v}_k \rangle + \frac{L \tilde{\theta}_k^2 }{2} \nonumber \\
			& = (1-\delta_k)\big[ f(\eta_k\mathbf{x}_k) - \Phi_k(\| \mathbf{x}^* \| \mathbf{v}_{k-1}) \big]  + \frac{L \delta_k^2 \| \mathbf{x}^* \|^2 }{2} \nonumber \\
			& = (1-\delta_k)\big[ f(\mathbf{x}_k) - \Phi_k(\| \mathbf{x}^* \| \mathbf{v}_{k-1}) \big]  + \frac{L \delta_k^2 \| \mathbf{x}^* \|^2 }{2} \nonumber
		\end{align}
		where the last inequality is because $f(\eta_k\mathbf{x}_k)\leq f(\mathbf{x}_k)$. For notational convenience, let $\mathbf{z}_{k+1}:= f(\mathbf{x}_{k+1}) - \Phi_{k+1}( \| \mathbf{x}^* \| \mathbf{v}_k)$, then the recursion is obvious, i.e., 
		\begin{align}\label{eq.apdx.recur}
			\mathbf{z}_{k+1} \leq (1-\delta_k) \mathbf{z}_k + 	\frac{L \delta_k^2 \| \mathbf{x}^* \|^2 }{2}.
		\end{align}
		This theorem can be proved by unrolling \eqref{eq.apdx.recur} and plugging in the choice of $\delta_k$.

\subsection{Proof of Lemma \ref{lemma.smooth_implification}}\label{apdx.lemma.smooth_implification}

Since $f$ is strictly convex and differentiable, we can define the Bregman divergence w.r.t. $f$ as
	\begin{align}
		D_f(\mathbf{x}, \mathbf{y}) = f(\mathbf{x}) - f(\mathbf{y}) -  \langle \nabla f(\mathbf{y}), \mathbf{x} - \mathbf{y} \rangle .
	\end{align}
	If we can show that $D_f(\mathbf{x}, \mathbf{y})\geq \frac{1}{2L} \|\nabla f(\mathbf{y}) - \nabla f(\mathbf{x})  \|_*^2$, this lemma can be proved.
	
	The strict convexity of $f$ implies that its conjugate $f^*$ that is also strictly convex.
	Using the duality of Bregman divergence, we have
	\begin{align*}
		D_f(\mathbf{x}, \mathbf{y}) = D_{f^*}\big(\nabla f(\mathbf{x}), \nabla f(\mathbf{y})\big). 
	\end{align*}
	Next by the fact that $f$ is $L$-smooth w.r.t. $\| \cdot \|$, it can be seen that $f^*$ is $\frac{1}{L}$-strongly convex w.r.t. $\| \cdot \|_*$ \cite{kakade2009duality}, from which one can conclude that
	\begin{align*}
		D_f(\mathbf{x}, \mathbf{y}) = D_{f^*}\big(\nabla f(\mathbf{x}), \nabla f(\mathbf{y})\big) \geq \frac{1}{2L} \|\nabla f(\mathbf{x})- \nabla f(\mathbf{y}) \|_*.
	\end{align*}
	The proof is thus completed.

\subsection{Proof of Theorem \ref {thm.dual}}\label{apdx.thm.dual}

Applying Lemma \ref{lemma.smooth_implification}, we have
	\begin{align*}
		& ~~~~~ \frac{1}{2L} \|\nabla f( \eta_k \mathbf{x}_k) - \nabla f(\mathbf{x}^*)  \|_*^2  \\
		& \leq f( \eta_k\mathbf{x}_k)  - f(\mathbf{x}^*) - \langle \nabla f(\mathbf{x}^*),  \eta_k\mathbf{x}_k - \mathbf{x}^* \rangle \\
		& \stackrel{(a)}{\leq} f( \eta_k\mathbf{x}_k) - f(\mathbf{x}^*) \leq \frac{2L\| \mathbf{x}^*\|^2}{k+1}
	\end{align*}
	where (a) is because of the optimality condition, i.e., $ \langle \nabla f(\mathbf{x}^*), \mathbf{x} - \mathbf{x}^* \rangle \geq 0, \forall \mathbf{x} \in {\cal K}$. This proves the first part of this theorem. 

For the second part, since we have $\nabla f(\mathbf{x}^*) \in {\cal K}^*$ (this can be obtained from the KKT condition), it follows directly that 
\begin{align*}
	\text{dist}^* \big(\nabla f(\mathbf{x}_\epsilon^*), {\cal K}^* \big)  \leq 	\|\nabla f( \eta_k \mathbf{x}_k) - \nabla f(\mathbf{x}^*)  \|_*.
\end{align*}
Hence the second part is proved.

\subsection{Proof of Theorem \ref{thm.stopcr}}\label{apdx.thm.stopcr}

By writing $\mathbf{g}_k$ as the summation of $\{ \nabla f(\eta_\tau \mathbf{x}_\tau)\}$, we have that
	\begin{align*}
		&~~~~ \| \mathbf{g}_k - \nabla f(\mathbf{x}^*) \|_* \\
		& \leq \sum_{\tau=0}^k \frac{2(\tau+1)}{(k+1)(k+2)} \| \nabla f( \eta_\tau \mathbf{x}_\tau) - \nabla f(\mathbf{x}^*) \|_*  \\
		& \leq 4 L \| \mathbf{x}^* \| \sum_{\tau=0}^k \frac{\sqrt{\tau+1}}{(k+1)(k+2)}  \\
		& \leq \frac{8 L \| \mathbf{x}^* \| \sqrt{k+2}}{3(k+1)}.
	\end{align*}
Then Using the fact that
	\begin{align*}
		\sqrt{k+2} \leq \sqrt{k+1} + \frac{1}{2\sqrt{k+1}}	
	\end{align*}
	we have
	\begin{align*}
		\| \mathbf{g}_k - \nabla f(\mathbf{x}^*) \|_* \leq \frac{8 L \| \mathbf{x}^* \| }{3\sqrt{k+1}} + \frac{4 L \| \mathbf{x}^* \| }{3(k+1)^{1.5}}.
	\end{align*}
	Considering the case where $k \geq 2$, we arrive at
	\begin{align*}
		&~~~~ \| \mathbf{g}_k - \nabla f(\mathbf{x}^*) \|_* \\
		& \leq \frac{8 L \| \mathbf{x}^* \| }{3\sqrt{k+1}} + \frac{4 L \| \mathbf{x}^* \| }{9\sqrt{k+1}} = \frac{28 L \| \mathbf{x}^* \|}{9\sqrt{k+1}}.
	\end{align*}
	Since we mush have $\nabla f(\mathbf{x}^*) \in {\cal K}^*$, the proof can be completed by simply taking square on both sides of the inequality above.

\subsection{Proof of Theorem \ref{thm.primal'}}\label{apdx.thm.primal'}
The proof of primal convergence in Theorem \ref{thm.primal} does not require strict convexity, and hence the convergence of $\eta_k \mathbf{X}_k$, i.e., Option II in Alg. \ref{alg.mcd_mem}, still holds
	\begin{align*}
		f \circ {\cal G} (\eta_{k+1} \mathbf{X}_{k+1}) -  f \circ {\cal G}(\mathbf{X}^*) \leq \frac{2L_{f \circ {\cal G}} \|\mathbf{X}^* \|^2}{k+2}.
	\end{align*}
	Comparing Options I and II, there is an underlying $\mathbf{X}_k$ associated with $\mathbf{y}_k$ satisfying $\mathbf{y}_k = {\cal G}(\mathbf{X}_k) -\mathbf{z}$ . Noticing that $f \circ {\cal G} (\eta_{k+1} \mathbf{X}_{k+1}) = f(\eta_{k+1}\mathbf{y}_k)$, the primal convergence is thus proved.
	
	The second inequality is established following a similar argument of \cite[Theorem C.1]{yurtsever2021}.

\subsection{Proof of Theorem \ref{thm.dual'}}\label{apdx.thm.dual'}
Let $\mathbf{y}^* := {\cal G}(\mathbf{X}^*) - \mathbf{z}$.
	Since $f$ is strictly convex and $L_f$ smooth, Lemma \ref{lemma.smooth_implification} still holds. Therefore, we have
	\begin{align}\label{eq.apdx.aaa}
		& ~~~~~ \frac{1}{2L_f} \|\nabla f( \eta_k \mathbf{y}_k) - \nabla f(\mathbf{y}^*)  \|_*^2  \\
		& \leq f( \eta_k\mathbf{y}_k)  - f(\mathbf{y}^*) - \langle \nabla f(\mathbf{y}^*),  \eta_k\mathbf{y}_k - \mathbf{y}^* \rangle \nonumber\\
		& \stackrel{(a)}{\leq} f( \eta_k\mathbf{y}_k) - f(\mathbf{y}^*) \nonumber\\
		& = f \circ {\cal G} (\eta_k \mathbf{X}_k) -  f \circ {\cal G}(\mathbf{X}^*) \nonumber \\
		& \leq \frac{2L_{f \circ {\cal G}}\| \mathbf{X}^*\|^2}{k+1} \nonumber
	\end{align}
	where (a) is because of the inequality below.
	\begin{align*}
		&~~\langle \nabla f(\mathbf{y}^*),  \eta_k\mathbf{y}_k - \mathbf{y}^* \rangle  \nonumber \\
		=& \langle 	\nabla f(\mathbf{y}^*), {\cal G}(\eta_k \mathbf{X}_k)- {\cal G}(\mathbf{X}^*) \rangle \\
		= & \langle 	{\cal G}^* \big(\nabla f(\mathbf{y}^*) \big), \eta_k \mathbf{X}_k- \mathbf{X}^*\rangle \geq 0
	\end{align*}
	where the first equation follows from $\langle {\cal G}^*(\mathbf{a}), \mathbf{X} \rangle = \langle \mathbf{a}, {\cal G}(\mathbf{X}) \rangle$; and the last inequality is because of the optimality condition of $\mathbf{X}^*$.
	
	Using Assumption \ref{as.4}, we have 
	\begin{align}\label{eq.apdx.aaaa}
		& ~~~~\| \nabla f \circ {\cal G}(\mathbf{X}_k) -  \nabla f \circ {\cal G}(\mathbf{X}^*) \|_* \\ 
		& = \| {\cal G}^* \big(\nabla f (\mathbf{y}_k) - \nabla f (\mathbf{y}^*) \big) \|_* \nonumber \\
		& \leq \bar{G} \| \nabla f (\mathbf{y}_k) - \nabla f (\mathbf{y}^*) \|_*. \nonumber
	\end{align}
	Combining \eqref{eq.apdx.aaa} and \eqref{eq.apdx.aaaa} completes the proof of first part of this theorem.

	For the second part, we have that
	\begin{align*}
		& ~~~~ \| \nabla f \circ {\cal G}(\hat{\mathbf{X}}_k) -  \nabla f \circ {\cal G}(\mathbf{X}^*) \|_* \nonumber \\
		& \leq \| \nabla f \circ {\cal G}(\hat{\mathbf{X}}_k) -  \nabla f \circ {\cal G}(\eta_k \mathbf{X}_k) \|_*  \\
		& ~~~~~~~~~~~~~~~~~~~~~~~~ + \| \nabla f \circ {\cal G}(\eta_k \mathbf{X}_k) -  \nabla f \circ {\cal G}(\mathbf{X}^*) \|_* \\
		& \leq L_{f \circ {\cal G}} \underbrace{\| \hat{\mathbf{X}}_k -  \eta_k \mathbf{X}_k \|}_{\text{first term}}  + \underbrace{\| \nabla f \circ {\cal G}(\eta_k \mathbf{X}_k) -  \nabla f \circ {\cal G}(\mathbf{X}^*) \|_*}_{\text{second term}}.
	\end{align*}
	The first term can be bounded using the same argument of \cite[Theorem C.1]{yurtsever2021}; see also Theorem \ref{thm.primal'}, where we have
	\begin{align*}
		\limsup_{k \rightarrow \infty } 	 \| \hat{\mathbf{X}}_k -  \eta_k \mathbf{X}_k \| \leq \limsup_{k  \rightarrow \infty} \mathbb{E}[\text{dist}(\hat{\mathbf{X}}_k, \Psi^*) ].
	\end{align*}
	The second term goes to $0$ when $k \rightarrow \infty $ following the first part of this proof.

\end{document}